\documentclass[a4,11pt]{amsart}
\usepackage{amssymb,amsmath,amsthm,txfonts}
\usepackage{cite}
\usepackage{color}
\usepackage[normalem]{ulem}

\newtheorem{thm}{Theorem}[section]
\newtheorem{lem}[thm]{Lemma}
\newtheorem{cor}[thm]{Corollary}
\newtheorem{prop}[thm]{Proposition}
\theoremstyle{definition}

\theoremstyle{remark}

\newtheorem{rems}[thm]{\textbf{Remarks}}

 \makeatletter
    
    \@addtoreset{equation}{section}
  \makeatother

\makeatletter
      \def\@makefnmark{%
         \leavevmode
            \raise.9ex\hbox{\check@mathfonts
                \fontsize\sf@size\z@\normalfont%
                            \@thefnmark}%
       }
      \makeatother

\newcommand{\p}{\mathbb{P}}

\newcommand{\D}{\textrm{div}}

\newcommand{\dd}{\textrm{d}}

\begin{document}

\title[]{Axisymmetric flows in the exterior of a cylinder}
\author[]{K. Abe}
\author[]{G. Seregin}
\date{}
\address[K. Abe]{Department of Mathematics, Graduate School of Science, Osaka City University, Sugimoto 3-3-138, Sumiyoshi-ku Osaka 558-8585, Japan}
\email{kabe@sci.osaka-cu.ac.jp}

\address[G. Seregin]{Mathematical Institute, Oxford University, Oxford, 24-29 St Giles’ OX1 3LB UK and Voronezh State University, Voronezh, Russia}
\email{seregin@maths.ox.ac.uk}

\date{\today}

\maketitle

\begin{abstract}
We study an initial-boundary value problem of the three-dimensional Navier-Stokes equations in the exterior of a cylinder $\Pi=\{x=(x_{h}, x_3)\ |\  |x_{h} |>1\}$, subject to the slip boundary condition. We construct unique global solutions for axisymmetric initial data $u_0\in L^{3}\cap L^{2}(\Pi)$ satisfying the decay condition of the swirl component $ru^{\theta}_{0}\in L^{\infty}(\Pi)$. 
\end{abstract}

\vspace{15pt}

%セクション1
\section{Introduction}
\vspace{10pt}
We consider the three-dimensional Navier-Stokes equations:

\begin{equation*}
\begin{aligned}
\partial_t u-\Delta{u}+u\cdot \nabla u+\nabla{p}&= 0    \\
\D\ u&=0 
\end{aligned}
\qquad \textrm{in}\ \Pi\times(0,\infty).
\tag{1.1}
\end{equation*}\\
It is well known that for small initial data $u_0\in L^{3}_{\sigma}(\mathbb{R}^{3})$, there exists a unique global solution $u\in BC([0,\infty); L^{3})$ of (1.1) \cite{Kato84}. However, unique existence of a global solution is unknown in general for large initial data in $L^{3}$ with finite energy. Here, $BC([0,\infty); X)$ denotes the space of all bounded and continuous functions from $[0,\infty)$ to a Banach space $X$ and $L^{p}_{\sigma}(\Pi)$ denotes the $L^{p}$-closure of compactly supported smooth solenoidal vector fields in a domain $\Pi\subset \mathbb{R}^{3}$. 

For initial data with finite energy $u_0\in L^{2}(\mathbb{R}^{3})$, it is well known that global Leray-Hopf weak solutions exist \cite{Leray1934}, \cite{Hopf51}. However, their regularity and uniqueness are unknown. For large initial data in $L^{3}(\mathbb{R}^{3})$, weak solutions are constructed in \cite{Calderon}, \cite{Lemarie}. See \cite{SS17} for weak $L^{3}$-solutions. 

The purpose of this paper is to construct unique global solutions of (1.1) for large axisymmetric initial data in $L^{3}\cap L^{2}$. We say that a vector field $u$ is axisymmetric if 

\begin{align*}
u(x)={}^{t}\hspace{-1pt}Ru(Rx)\quad x\in \mathbb{R}^{3},\ \eta\in [0,2\pi],
\end{align*}\\
for $R=(e_{r}(\eta),e_{\theta}(\eta),e_{z})$ and $e_{r}(\eta)={}^{t}(\cos\eta,\sin\eta,0)$, $e_{\theta}(\eta)={}^{t}(-\sin\eta,\cos\eta,0)$, $e_{z}={}^{t}(0,0,1)$. We say that a scaler function $p$ is axisymmetric if $p(x)=p(Rx)$ for $x\in \mathbb{R}^{3}$ and $\eta\in [0,2\pi]$. We set the cylindrical coordinate $(r,\theta,z)$ by $x_1=r\cos\theta$, $x_2=r\sin\theta$, $x_{3}=z$ and decompose the axisymmetric vector field into three terms:

\begin{align*}
u(x)=u^{r}(r,z)e_{r}(\theta)+u^{\theta}(r,z)e_{\theta}(\theta)+u^{z}(r,z)e_{z}.
\end{align*}\\
The azimuthal component $u^{\theta}$ is called swirl velocity (see, e.g., \cite{MaB}).

Unique global solutions of (1.1) for axisymmetric initial data without swirl were first constructed in \cite{La68b}, \cite{UI} by the Galerkin approximation. Later on, unique global solutions are constructed in \cite{LMNP} by a strong solution approach for axisymmetric data without swirl in $H^{2}(\mathbb{R}^{3})$. See also \cite{Abidi08} for $H^{1/2}(\mathbb{R}^{3})$. For axisymmetric solutions of (1.1), the vorticity $\omega=\textrm{curl}\ u$ is expressed by 

\begin{align*}
\omega
&=\omega^{r}e_{r}+\omega^{\theta}e_{\theta}+\omega^{z}e_{z}\\
&=(-\partial_z u^{\theta})e_{r}+(\partial_{z}u^{r}-\partial_{r}u^{z})e_{\theta}+\Big(\partial_ru^{\theta}+\frac{u^{\theta}}{r}\Big)e_{z},
\end{align*}\\
and for $v=u^{r}e_{r}+u^{z}e_{z}$, the azimuthal component $\omega^{\theta}$ satisfies the vorticity equation

\begin{equation*}
\begin{aligned}
\partial_t \Big(\frac{\omega^{\theta}}{r}\Big)+v\cdot \nabla \Big(\frac{\omega^{\theta}}{r}\Big)-\Big(\Delta +\frac{2}{r}\partial_r \Big)\Big(\frac{\omega^{\theta}}{r}\Big) &=\partial_z\Big(\frac{  u^{\theta}}{r}\Big)^{2}.
\end{aligned}
\tag{1.2}
\end{equation*}\\
For axisymmetric solutions without swirl, the right-hand side vanishes and the global a priori estimate 

\begin{align*}
\Big\|\frac{\omega^{\theta}}{r}\Big\|_{L^{2}(\mathbb{R}^{3})}\leq \Big\|\frac{\omega^{\theta}_{0}}{r}\Big\|_{L^{2}(\mathbb{R}^{3})}\quad t>0,
\end{align*}\\
holds. The above vorticity estimate implies existence of unique global solutions for axisymmetric data without swirl $u_0\in L^{3}\cap L^{2}(\mathbb{R}^{3})$. (We may assume the condition $\omega^{\theta}_{0}/r\in L^{2}(\mathbb{R}^{3})$ since local-in-time solutions belong to $H^{2}(\mathbb{R}^{3})$.) In other words, unique global solutions exist for large axisymmetric initial data in $L^{3}\cap L^{2}(\mathbb{R}^{3})$, provided that without swirl. For axisymmetric data with swirl, unique existence of global solutions in $\mathbb{R}^{3}$ is unknown.

In this paper, we study axisymmetric solutions \textit{with swirl} in the exterior of a cylinder

\begin{align*}
\Pi=\{x=(x_{1}, x_2,x_3 )\in \mathbb{R}^{3}\ |\ |x_{h}|>1,\ x_{h}=(x_1,x_2)   \},
\end{align*}\\
subject to the slip boundary condition

\begin{align*}
(D(u)n)_{\textrm{tan}}=0,\ u\cdot n=0\quad \textrm{on}\quad \partial\Pi.    \tag{1.3}
\end{align*}\\
Here, $n=-e_{r}$ denotes the unit outward normal vector field on $\partial\Pi$, $D(u)=(\nabla u+\nabla^{T}u)/2$ is the deformation tensor and $f_{\textrm{tan}}=f-n(f\cdot n)$ is a tangential component of a vector field $f$ on $\partial\Pi$. Since axisymmetric vector fields $u=u^{r}e_{r}+u^{\theta}e_{\theta}+u^{z}e_{z}$ satisfy

\begin{align*}
u^{r}=0,\quad \partial_r u^{\theta}-u^{\theta}=0,\quad \partial_r u^{z}=0\quad \textrm{on}\ \{r=1\},   
\end{align*}\\
subject to the slip boundary condition (1.3), the azimuthal component of vorticity $\omega^{\theta}$ vanishes on the boundary (see Remarks 6.1 (ii) for the Dirichlet boundary condition). 

By the partial regularity result \cite{CKN}, it is expected that axisymmetric solutions are smooth in the interior of $\Pi$. Moreover, as noted in \cite{CHKLSY}, they will not develop singularities on the boundary due to viscosity. See \cite{Seregin02}, \cite{GKT06} for partial regularity results up to the boundary subject to the Dirichlet boundary condition. The regularity theory for the slip boundary condition (1.3) may be simpler than that for the Dirichlet boundary condition. In fact, for a half space regularity results are deduced from a whole space case by a reflection argument; see \cite{BaeJin08}. In this paper, we prove that axisymmetric solutions are sufficiently smooth in the exterior of the cylinder $\overline{\Pi}\times (0,\infty)$, subject to the slip boundary condition (1.3). We impose the slip boundary condition in order to construct approximate solutions for $\mathbb{R}^{3}$; see Remarks 1.2 (iii).

Our goal is to construct unique global mild solutions of (1.1) for axisymmetric initial data with swirl in $L^{3}\cap L^{2}(\Pi)$. Since the boundary of the cylinder $\Pi\subset \mathbb{R}^{3}$ is uniformly regular, we construct mild solutions by using the $\tilde{L}^{p}$-theory. We set 

\begin{align*}
\tilde{L}^{p}(\Pi)=L^{p}\cap L^{2}(\Pi)
\end{align*}\\
(resp. $\tilde{L}^{p}_{\sigma}(\Pi)=L^{p}_{\sigma}\cap L^{2}_{\sigma}(\Pi)$) for $p\in [2,\infty)$. It is proved in \cite{FKS1} (\cite{FKS2}) that that the Helmholtz projection $\mathbb{P}$ acts as a bounded operator on $\tilde{L}^{p}(\Pi)$. Moreover, it is recently shown in \cite{FR} that the Stokes operator subject to the slip boundary condition $A=\mathbb{P}\Delta$ generates a $C_{0}$-analytic semigroup on $\tilde{L}^{p}_{\sigma}(\Pi)$ (see also \cite{FKS1}, \cite{FKS3} for the Dirichlet boundary condition). We construct mild solutions for $u_0\in \tilde{L}^{3}_{\sigma}(\Pi)$ of the form

\begin{align*}
u(t)=e^{tA}u_0-\int_{0}^{t}e^{(t-s)A}\p (u\cdot \nabla u)(s)\dd s.   \tag{1.4}
\end{align*} \\
Since the swirl component satisfies the Robin boundary condition, axisymmetric solutions of (1.4) satisfy the energy equality

\begin{align*}
\int_{\Pi}|u|^{2}\dd x+2\int_{0}^{t}\hspace{-3pt}\int_{\Pi}\Big(|\nabla v|^{2}+|\nabla u^{\theta}|^{2}+\Big|\frac{u^{\theta}}{r}\Big|^{2}\Big)\dd x\dd s
+2 \int_{0}^{t}\hspace{-3pt}\int_{\partial \Pi}|u^{\theta}|^{2}\dd {\mathcal{H}} \dd s= \int_{\Pi}|u_{0}|^{2}\dd x,   \tag{1.5}
\end{align*}\\
where $\dd {\mathcal{H}}$ denotes the surface element on $\partial\Pi$.

We construct unique global solutions for large axisymmetric data with swirl $u_0\in \tilde{L}^{3}_{\sigma}(\Pi)$ satisfying the decay condition of the swirl component $ru^{\theta}_{0}\in L^{\infty}(\Pi)$. The main result of this paper is the following:

\vspace{15pt}

%thm1 
\begin{thm}
Let $u_0\in \tilde{L}^{3}_{\sigma}(\Pi)$ be an axisymmetric vector field. Assume that $ru^{\theta}_{0}\in L^{\infty}(\Pi)$. Then, there exists a unique axisymmetric mild solution $u\in BC([0,\infty); \tilde{L}^{3}(\Pi))$ satisfying (1.5) for $t\geq 0$.
\end{thm}

\vspace{10pt}

\begin{rems}
(i) It is unknown in general whether axisymmetric solutions in $\mathbb{R}^{3}$ for $u_0\in \tilde{L}^{3}_{\sigma}(\mathbb{R}^{3})$ satisfying $ru^{\theta}_{0}\in L^{\infty}(\mathbb{R}^{3})$ are globally bounded for all $t>0$. See \cite{NP1}, \cite{NP2}, \cite{CL}, \cite{JX} for regularity criteria of axisymmetric solutions. For axisymmetric solutions, an upper bound of the form $|u(x,t)|\leq Cr^{-1}$, $r<1$, is called type I condition. It is proved in \cite{CSTY1}, \cite{CSTY2} by De Giorgi method and \cite{KNSS}, \cite{SS} by the Liouville-type theorem that axisymmetric solutions do not develop type I singularities. See  \cite{Se15} about type I singularities. Recently, it is shown in \cite{LNZ} (\cite{LZ}) that axisymmetric smooth solutions in $\mathbb{R}^{3}\times(-T,0)$ for $u(\cdot,-T)\in L^{2}(\mathbb{R}^{3})$ and $ru^{\theta}(\cdot,-T)\in L^{\infty}(\mathbb{R}^{3})$ satisfy an upper bound of the form $|u(x,t)|\leq C|\log {r}|^{1/2}r^{-2}$ near $(r,t)=0$ with some constant $C$.

\noindent
(ii) It is known that solutions of (1.1) in $\mathbb{R}^{3}$ are smooth if the direction of vorticity is Lipschitz continuous for spatial variables in regions of high vorticity magnitude \cite{CF} (called a geometric regularity criterion). For axisymmetric flows without swirl, vorticity varies only in the azimuthal direction and is identified with a scalar function. On the other hand, for axisymmetric flows with swirl vorticity varies also in the radial and vertical directions. We constructed unique global solutions whose vorticity may become large and vary in three directions. For a half space $\mathbb{R}^{3}_{+}$, a geometric regularity criterion is proved in \cite{daVeiga06}, subject to the slip boundary condition. See also \cite{daVeiga07} for the Dirichlet boundary condition.

\noindent
(iii) Theorem 1.1 implies existence of approximate solutions for $\mathbb{R}^{3}$. Since the exterior of the cylinder $\Pi^{\varepsilon}=\{r>\varepsilon\}$ approaches $\mathbb{R}^{3}$ as $\varepsilon\to 0$, axisymmetric solutions in $\mathbb{R}^{3}$ can be viewed as limits of solutions in $\Pi^{\varepsilon}$. Indeed, axisymmetric solutions without swirl in $\Pi^{\varepsilon}$ are uniformly bounded in $L^{\infty}_{t}H^{1}_{x}$ for $\varepsilon>0$ and approach those in $\mathbb{R}^{3}$ \cite[p.78, l.7]{La68b}. See Remarks 6.1 (iii). For the case with swirl, unique existence of global solutions is proved in \cite{Za} (\cite{Zaja07}) in a bounded cylindrical domain for sufficiently smooth initial data. It is unknown whether global solutions with swirl are uniformly bounded for all $\varepsilon>0$. We constructed unique global mild solutions for $u_0\in \tilde{L}^{3}_{\sigma}(\Pi^{\varepsilon})$ satisfying the uniform estimate for the swirl component (1.6). 
\end{rems}

\vspace{20pt}

Let us sketch the proof of Theorem 1.1. We first construct local-in-time mild solutions of (1.4) for $u_0\in \tilde{L}^{3}_{\sigma}$ and prove that mild solutions are axisymmetric and satisfy the energy equality (1.5) for axisymmetric initial data. The major step of the proof is to derive a global $L^{4}$-bound for axisymmetric solutions $u=v+u^{\theta}e_{\theta}$. Once we obtain the global bound, it is not difficult to see that $u\in BC([0,\infty); \tilde{L}^{3})$ by local solvability and the energy equality (1.5).

We first prove the global $L^{\infty}$-estimate for the swirl component

\begin{align*}
||ru^{\theta}||_{L^{\infty}(\Pi)}\leq ||ru^{\theta}_{0}||_{L^{\infty}(\Pi)}\quad t>0.  \tag{1.6}
\end{align*}\\
Since $r\geq 1$ in the exterior of the cylinder $\Pi$, the $L^{\infty}$-estimate (1.6) and the energy equality (1.5) implies the global $L^{4}$-bound for $u^{\theta}$ of the form 

\begin{align*}
||u^{\theta}||_{4}\leq ||ru^{\theta}_{0}||_{\infty}^{\frac{1}{2}}||u_0||_{2}^{\frac{1}{2}}\quad t>0. \tag{1.7}
\end{align*}\\
In order to prove (1.6), we study the drift-diffusion equation subject to the Robin boundary condition:

\begin{equation*}
\begin{aligned}
\partial_t \Gamma +b\cdot \nabla \Gamma -\Delta\Gamma +\frac{2}{r}\partial_r \Gamma&=0\hspace{15pt} \textrm{in}\ \Pi\times (0,T), \\
\partial_{n} \Gamma+2\Gamma&=0\hspace{15pt} \textrm{on}\ \partial\Pi\times (0,T),\\
\Gamma&=\Gamma_0\hspace{10pt} \textrm{on}\ \Pi\times \{t=0\}.
\end{aligned}
\tag{1.8}
\end{equation*}\\
Here, $\partial_n=-\partial_r$ denotes the normal derivative. The function $\Gamma=ru^{\theta}$ is a solution of (1.8) for $b=v$. We prove the $L^{\infty}$-estimate

\begin{align*}
||\Gamma||_{L^{\infty}(\Pi)}\leq ||\Gamma_{0}||_{L^{\infty}(\Pi)}\quad t>0\textcolor{red}{,} \tag{1.9}
\end{align*}\\
for solutions to (1.8). Since the sign of the coefficient is plus in the Robin boundary condition, a maximum principle holds if the coefficient $b$ and $\Gamma$ are bounded in $\Pi\times[0,T]$. Then the $L^{\infty}$-estimate (1.9) easily follows from a maximum principle (see Lemma 3.1). If $\Gamma$ is decaying sufficiently fast as $|x|\to\infty$, we are able to obtain (1.9) by estimating $L^{p}$-norms of $\Gamma$ for $p=2^{m}$ and sending $m\to\infty$. Since we assume that $ru^{\theta}_0$ is merely bounded, the function $ru^{\theta}$ may not decay as $|x|\to\infty$ . We shall prove (1.9) for non-decaying solutions $\Gamma$.

We apply the $L^{\infty}$-estimate (1.9) for $ru^{\theta}$ and obtain (1.6). Note that the boundedness of $ru^{\theta}$ does not follow from properties of local-in-time solutions to (1.1) for $u_0\in \tilde{L}^{3}_{\sigma}$. For this purpose, we first extend the $L^{\infty}$-estimate (1.9) for mild solutions to (1.8) for $\Gamma_0\in L^{\infty}$ and the coefficient $b$ such that $t^{1/2-3/2p}b\in C([0,T]; L^{p})$ vanishes at time zero for $p\in (3,\infty]$. We then deduce from the integral form (1.4) that $ru^{\theta}$ is a mild solution to (1.8) (see Lemma 4.7). 

We next estimate a global $L^{4}$-norm of $v=u^{r}e_{r}+u^{z}e_z$. We apply an interpolation inequality 

\begin{align*}
||v||_{4}\leq C||v||_{2}^{\frac{1}{4}}(||v||_{2}+||\omega^{\theta}||_{2})^{\frac{3}{4}},   \tag{1.10}
\end{align*}\\
and estimate an energy norm of the vorticity $\omega^{\theta}$. Since $\omega^{\theta}$ vanishes on the boundary, we control the external force $\partial_z(u^{\theta}/r)^{2}$ by using viscosity and estimate 

\begin{equation*}
\begin{aligned}
\int_{\Pi}\Big|\frac{\omega^{\theta}}{r}\Big|^{2}\dd x
+\int_{0}^{t}\int_{\Pi}\Big|\nabla \Big(\frac{\omega^{\theta}}{r}\Big)\Big|^{2}\dd x\dd s
&\leq \int_{\Pi}\Big|\frac{\omega^{\theta}_{0}}{r}\Big|^{2}\dd x
+||ru_{0}^{\theta}||_{\infty}^{2}||u_0||_{2}^{2}\\
&=: E\qquad t>0.
\end{aligned}
\tag{1.11}
\end{equation*}\\
Since the above vorticity estimate implies the global bound

\begin{equation*}
\begin{aligned}
\int_{\Pi}|\omega^{\theta}|^{2}\dd x
+\int_{0}^{t}\int_{\Pi}\Bigg(|\nabla\omega^{\theta}|^{2}+\Big|\frac{\omega^{\theta}}{r}\Big|^{2}   \Bigg)\dd x\dd s
&\leq \int_{\Pi}|\omega^{\theta}_{0}|^{2}\dd x  \\
&+C(E^{\frac{3}{4}}||u_0||_{2}^{\frac{1}{2}}
+||ru^{\theta}_{0}||_{\infty}^{2} )||u_0||_{2}^{2},\quad t>0,
\end{aligned}
\tag{1.12}
\end{equation*}\\
the local-in-time solution $u=v+u^{\theta}e_{\theta}$ is globally bounded on $L^{4}$.

\vspace{15pt}

This paper is organized as follows. In Section 2, we state a local existence theorem of mild solutions for $u_0\in \tilde{L}^{3}_{\sigma}$ and prove axial symmetry of mild solutions. In Section 3, we study the drift-diffusion equation (1.8) for a bounded coefficient and prove the $L^{\infty}$-estimate (1.9) by a maximum principle. In Section 4, we extend (1.9) for mild solutions to (1.8) under the weak regularity condition of a coefficient, and apply (1.9) for the swirl component of axisymmetric solutions. In Section 5, we prove the a priori estimates (1.11) and (1.12). In Section 6, we prove Theorem 1.1. In Appendix A, we give a proof for a local solvability result stated in Section 2. In Appendix B, we prove some interpolation inequalities used in Section 5.

\vspace{20pt}

%セクション2
\section{Local existence of axisymmetric solutions on $\tilde{L}^{3}$}

\vspace{10pt}

In this section, we construct local-in-time axisymmetric solutions of (1.1) for $u_0\in \tilde{L}^{3}_{\sigma}$ satisfying the energy equality (1.5). Local solvability for $u_0\in \tilde{L}^{3}_{\sigma}$ is known for $\mathbb{R}^{3}$ \cite[Theorem 3]{Kato84}. We give a proof for the exterior of the cylinder by using $\tilde{L}^{p}$-theory in Appendix A.

\vspace{15pt}

\subsection{Local solvability}

Let $C^{\alpha}([\delta,T]; X)$ denote the space of all $\alpha$-th H\"older continuous functions $f\in C([\delta, T]; X )$ for a Banach space $X$. Let $C^{\alpha}((0,T]; X)$ denote the space of functions in $C^{\alpha}([\delta,T]; X)$ for all $\delta\in (0,T)$. For the convenience, we denote by $\tilde{L}^{p}=L^{p}\cap L^{2}$ also for $p=\infty$.

\vspace{15pt}

%2.1
\begin{lem}
For $u_0\in \tilde{L}^{3}_{\sigma}$, there exist $T>0$ and a unique mild solution of (1.4) satisfying  

\begin{align*}
&t^{\frac{3}{2}(\frac{1}{3}-\frac{1}{p})}u\in C([0,T]; \tilde{L}^{p}),\quad  3\leq p\leq \infty,  \tag{2.1}   \\
&t^{\frac{3}{2}(\frac{1}{3}-\frac{1}{r})+\frac{1}{2}}\nabla u\in C([0,T]; \tilde{L}^{r}),\quad  3\leq r< \infty,    \tag{2.2}
\end{align*}\\
$t^{3/2(1/3-1/p)}u$ and $t^{3/2(1/3-1/r)+1/2}\nabla u$ vanish at time zero except for $p=3$. Moreover, 

\begin{equation*}
\begin{aligned}
&u\in C^{\alpha}((0,T]; \tilde{L}^{3}),\\
&\nabla u\in C^{\frac{\alpha}{2}}((0,T]; \tilde{L}^{3}),\quad 0<\alpha<1.\\
\end{aligned}
\tag{2.3}
\end{equation*}
\end{lem}

\vspace{15pt}
We show that mild solutions satisfy (1.1) by applying an abstract regularity result \cite[4.3.1 Theorem 4.3.4]{Lunardi}.

\vspace{15pt}

%2.2
\begin{prop}
Let $B$ be a generator of an analytic semigroup in a Banach space $X$ with a domain $D(B)$. Assume that $f\in L^{1}(0,T; X)\cap C^{\beta}((0,T]; X)$ for $\beta\in (0,1)$. Then, 

\begin{align*}
w=\int_{0}^{t}e^{(t-s)B}f(s)\dd s
\end{align*}\\
belongs to $C^{\beta}((0,T]; D(B))\cap C^{1+\beta}((0,T]; X)$.
\end{prop}

\vspace{5pt}

%2.3
\begin{prop}
The mild solution $u$ in Lemma 2.1 satisfies 

\begin{align*}
u\in C^{\gamma}((0,T]; D(A))\cap C^{1+\gamma}((0,T]; L^{2}),\quad 0<\gamma<\frac{1}{2},  \tag{2.4}
\end{align*}\\
for $D(A)=\{u\in L^{2}_{\sigma}\cap H^{2}\ |\ (D(u)n)_{\textrm{tan}}=0,\ u\cdot n=0\ \partial\Pi\ \}$. In particular, $u$ satisfies the equations (1.1) and (1.3).
\end{prop}

\vspace{5pt}

\begin{proof}
We set $f=-\mathbb{P}u\cdot \nabla u$. It follows from (2.1)-(2.3) that 

\begin{align*}
||f||_{2}
&\leq ||u\cdot \nabla u||_{2}
\leq ||u||_{3}||\nabla u||_{6}
\leq \frac{C}{t^{\frac{3}{4}}},\\
||f(t)-f(\tau)||_{2}
&\leq ||(u(t)-u(\tau))\cdot \nabla u(t)||_{2}
+||u(\tau)\cdot \nabla (u(t)-u(\tau))||_{2}\\
&\leq ||u(t)-u(\tau)||_{3}||\nabla u(t)||_{6}
+||u(\tau)||_{6}||\nabla u(t)-\nabla u(\tau)||_{3}\\
&\leq C\Big(\frac{|t-\tau|^{\alpha}}{t^{\frac{3}{4}}}+\frac{|t-\tau|^{\frac{\alpha}{2}}}{\tau^{\frac{1}{4}}}  \Big)\qquad \textrm{for}\quad  0<\tau<t\leq T.
\end{align*}\\
Thus $f\in L^{1}(0,T; L^{2})\cap C^{\alpha/2}((0,T]; L^{2})$ for $\alpha\in (0,1)$. Applying Proposition 2.2 yields (2.4).
\end{proof}

\vspace{5pt}

%subsection 2.2
\subsection{Axial symmetry}

We show that mild solutions are axisymmetric and satisfies the energy equality (1.5) for axisymmetric initial data.

\vspace{10pt}

%lem2.4
\begin{lem}
Assume that $u_0$ is axisymmetric. Then, the mild solution $u$ in Lemma 2.1 is axisymmetric and satisfies 

\begin{equation*}
\begin{aligned}
\partial_t u^{r}+v\cdot \nabla u^{r}-\frac{|u^{\theta}|^{2}}{r}
-\Big(\Delta-\frac{1}{r^{2}}\Big)u^{r}+\partial_r p&=0\\
\partial_t u^{\theta}+v\cdot \nabla u^{\theta}+\frac{u^{r}}{r}u^{\theta}
-\Big(\Delta-\frac{1}{r^{2}}\Big)u^{\theta}&=0\\
\partial_t u^{z}+v\cdot \nabla u^{z}-\Delta u^{z}+\partial_z p&=0\\
\partial_r u^{r}+\frac{u^{r}}{r}+\partial_z u^{z}&=0
\end{aligned}
\qquad \textrm{in}\ \Pi\times (0,T), \tag{2.5}
\end{equation*}
\begin{align*}
u^{r}=0,\quad \partial_r u^{\theta}-u^{\theta}=0,\quad \partial_r u^{z}=0\quad \textrm{on}\ \partial\Pi\times (0,T),   \tag{2.6}
\end{align*}\\
and the energy equality (1.5).
\end{lem}

\vspace{5pt}

%prop2.3
\begin{prop}
Assume that a vector field $u=u^{r}e_{r}+u^{\theta}e_{\theta}+u^{z}e_{z}$ satisfies (1.3). Then, $(u^{r},u^{\theta},u^{z})$ satisfies (2.6). The converse also holds. 
\end{prop}

\vspace{5pt}

\begin{proof}
By fundamental calculations using the cylindrical coordinate, we observe that 

\begin{align*}
D(u^{r}e_{r})e_{r}&=\partial_r u^{r}e_{r}+\frac{1}{2r}\partial_{\theta}u^{r}e_{\theta}+\frac{1}{2}\partial_z u^{r}e_{z},\\
D(u^{\theta}e_{\theta})e_{r}&=\frac{1}{2}\Big(\partial_r u^{\theta}-\frac{u^{\theta}}{r}\Big)e_{\theta},\\
D(u^{z}e_{z})e_{r}&=\frac{1}{2}\partial_r u^{z}e_{z},\\
D(u)e_{r}&=\partial_r u^{r}e_{r}+\frac{1}{2}\Big(\frac{1}{r}\partial_{\theta}u^{r}+\partial_ru^{\theta}-\frac{u^{\theta}}{r}\Big)e_{\theta}+\frac{1}{2}(\partial_z u^{r}+\partial_r u^{z})e_{z}.
\end{align*}\\
By (1.3), $(u^{r},u^{\theta},u^{z})$ satisfies (2.6). Conversely, suppose that (2.6) holds. Then, 

\begin{equation*}
D(u)e_{r}=\partial_ru^{r}e_{r},\quad u\cdot e_{r}=0  \qquad \textrm{on}\ \{r=1\}.
\end{equation*}\\
Thus (1.3) holds for $u=u^{r}e_{r}+u^{\theta}e_{\theta}+u^{z}e_{z}$.
\end{proof}

\vspace{15pt}

%prop 2.6
\begin{prop}
Set the rotation operator $U=U_{\eta}: L^{2}(\Pi)\longrightarrow L^{2}(\Pi)$ by 

\begin{align*}
f(x)\longmapsto {}^{t}Rf(Rx)
\end{align*}\\
and $R=(e_{r}(\eta),e_{\theta}(\eta),e_{z})$ for $\eta\in [0,2\pi]$. Then, we have 

\begin{align*}
Ue^{tA}f&=e^{tA}Uf,    \tag{2.7}\\
U\mathbb{P}g&=\mathbb{P}Ug,   \tag{2.8}\\
U(h\cdot \nabla h)&=(Uh)\cdot \nabla (Uh), \tag{2.9}
\end{align*}\\
for $f\in L^{2}_{\sigma}$, $g\in L^{2}$ and $h\in H^{1}$ satisfying $h\cdot \nabla h\in L^{2}$.
\end{prop}

\vspace{5pt}

\begin{proof}
We give a proof for (2.7). We are able to prove (2.8) and (2.9) by a similar way. We set $w=e^{tA}f$ and $w_{\eta}=U_{\eta}w$. Since the Stokes equations are rotationally invariant, $w_{\eta}$ satisfies 

\begin{equation*}
\begin{aligned}
\partial_t w_{\eta}-\Delta w_{\eta}+\nabla q_{\eta}&=0\\
\D\ w_{\eta}&=0
\end{aligned}
\quad \textrm{in}\ \Pi\times (0,\infty),
\end{equation*}\\
with some associated pressure $q_{\eta}$. It follows that 

\begin{align*}
w_{\eta}(x)
&={}^{t}Rw(Rx)\\
&= w^{r}(r,\theta+\eta,z)e_{r}(\theta)+w^{\theta}(r,\theta+\eta,z)e_{\theta}(\theta)+w^{z}(r,\theta+\eta,z)e_{z}.
\end{align*}\\
Since $(w^{r}, w^{\theta},w^{z})$ satisfies (2.6) by Proposition 2.5, $w_{\eta}$ satisfies the slip boundary condition (1.3). Since $w_{\eta}$ is a unique solution of the Stokes equations for $f_{\eta}=U_{\eta}f$, we have $w_{\eta}=e^{tA}f_{\eta}$. 
\end{proof}

\vspace{5pt}

\begin{proof}[Proof of Lemma 2.4]
We multiply $U$ by (1.4). It follows from (2.7)-(2.9) that 

\begin{align*}
Uu&=Ue^{tA}u_0-\int_{0}^{t}Ue^{(t-s)A}\mathbb{P}(u\cdot \nabla u)(s)\dd s\\
&=e^{tA}Uu_0-\int_{0}^{t}e^{(t-s)A}\mathbb{P}(Uu\cdot \nabla Uu)(s)\dd s.
\end{align*}\\
Since $u_0$ is axisymmetric, $u_0=Uu_0$. Hence $Uu$ is a mild solution of (1.1) for $u_0$. By the uniqueness of the mild solution, we have $u=U_{\eta}u$ for $\eta\in [0,2\pi]$. Thus $u$ is axisymmetric. Since $u$ satisfies (1.1) and (1.3) by Proposition 2.3, $(u^{r},u^{\theta},u^{z})$ satisfies (2.5) and (2.6). The energy equality (1.5) follows from integration by parts.
\end{proof}

\vspace{20pt}

%セクション3
\section{A maximum principle}

\vspace{10pt}
We consider the drift-diffusion equation (1.8) with a bounded coefficient and prove the $L^{\infty}$-estimate (1.9) by a maximum principle. Let $C(\overline{\Pi}\times [0,T])$ denote the space of all bounded and continuous functions in $\overline{\Pi}\times [0,T]$. Let $C^{2,1}(\overline{\Pi}\times [\delta,T])$ denote the space of all functions $f\in C(\overline{\Pi}\times [\delta,T])$ such that $\partial_t^{s}\partial_x^{k}f\in C(\overline{\Pi}\times [\delta,T])$ for $2s+|k|\leq 2$. We denote by $C^{2,1}(\overline{\Pi}\times (0,T])$ the space of all functions in $C^{2,1}(\overline{\Pi}\times [\delta ,T])$ for all $\delta \in (0,T)$. The goal of this section is: 

\vspace{15pt}

%lem3.1
\begin{lem}
Let $\Gamma\in C^{2,1}(\overline{\Pi}\times (0,T])\cap C(\overline{\Pi}\times [0,T])$ be a solution of (1.8). Assume that $b\in C(\overline{\Pi}\times [0,T])$. Then, the $L^{\infty}$-estimate (1.9) holds for $t\geq 0$.
\end{lem}

\vspace{15pt}

We prove Lemma 3.1 by a maximum principle. When $\Pi$ is bounded, a maximum principle with the Robin boundary condition is known \cite[Lemma 2.3]{Lieb}. We give a proof for the unbounded domain $\Pi$.

\vspace{15pt}

%prop3.2
\begin{prop}[Maximum principle]
Assume that $\Gamma\in C^{2,1}(\overline{\Pi}\times (0,T])\cap C(\overline{\Pi}\times [0,T])$ satisfies

\begin{align*}
\partial_t \Gamma +b\cdot \nabla \Gamma -\Delta \Gamma +\frac{2}{r}\partial_r \Gamma&\leq 0\quad \textrm{in}\ \Pi\times (0,T],   \tag{3.1}\\
\partial_n\Gamma +2\Gamma &\leq 0\quad \textrm{on}\ \partial \Pi\times (0,T], \tag{3.2}  \\
\Gamma&\leq 0\quad \textrm{on}\ \Pi\times \{t=0\}.   \tag{3.3} 
\end{align*}\\
Then, 

\begin{align*}
\Gamma\leq 0\quad \textrm{in}\ \Pi\times [0,T]. \tag{3.4}
\end{align*}
\end{prop}

\vspace{10pt}

%cor3.3
\begin{cor}
Assume that the reverse inequalities of (3.1)-(3.3) hold. Then, $\Gamma\geq 0$ in $\Pi\times [0,T]$.
\end{cor}

\vspace{5pt}

\begin{proof}[Proof of Lemma 3.1]
We set 

\begin{align*}
M&=\sup_{x\in \Pi}\Gamma_0(x),\\
m&=\inf_{x\in \Pi}\Gamma_0(x).
\end{align*}\\
We first show (1.9) when $m\leq 0$. We set 

\begin{align*}
\Gamma_m=m-\Gamma.
\end{align*}\\
The function $\Gamma_m$ satisfies (3.1) and (3.3). Since $m\leq 0$, it follows that 

\begin{align*}
(\partial_n+2)\Gamma_m
&=2m-(\partial_n+2)\Gamma\\
&=2m\leq 0.
\end{align*} \\
Hence the condition (3.2) is satisfied. Applying Proposition 3.2 implies that 

\begin{align*}
m\leq \Gamma(x,t)\quad \textrm{in}\ \Pi\times [0,T].   \tag{3.5}
\end{align*} \\
We next estimate $\Gamma$ from above. We first consider the case $M\leq 0$. Since $\Gamma_0\leq M \leq 0$, we apply Proposition 3.2 to $\Gamma$ and observe that $\Gamma\leq 0$. It follows from (3.5) that 

\begin{align*}
||\Gamma||_{\infty}
&=-\inf_{x\in \Pi}\Gamma (x,t)\\
&\leq -m=||\Gamma_0||_{\infty}.
\end{align*}\\
Thus (1.9) holds. We next consider the case $M>0$. We set 

\begin{align*}
\Gamma_M=M-\Gamma.
\end{align*}\\
Since $(\partial_n+2)\Gamma_M=2M>0$, the reverse inequalities of (3.1)-(3.3) hold for $\Gamma_M$. Applying Corollary 3.3 implies that 

\begin{align*}
\Gamma(x,t)\leq M\quad \textrm{in}\ \Pi\times [0,T].  \tag{3.6}
\end{align*}\\
By (3.5) and (3.6), we obtain  

\begin{align*}
||\Gamma||_{\infty}
&=\max\Big\{-\inf_{x\in \Pi} \Gamma(x,t),\ \sup_{x\in \Pi} \Gamma(x,t) \Big\}\\
&\leq \max\{-m,\ M  \}= ||\Gamma_0||_{\infty}.
\end{align*}\\
We proved (1.9) when $m\leq 0$.\\

It remains to show (1.9) when $m>0$. Since $\Gamma_0\geq m>0$, we observe that $\Gamma \geq 0$ by Corollary 3.3. Applying Corollary 3.3 for $\Gamma_{M}=M-\Gamma$ implies that $0\leq \Gamma\leq M$. Thus (1.9) holds when $m>0$. The proof is complete.
\end{proof}

\vspace{15pt}
We prove Proposition 3.2 from the following:
\vspace{15pt}

%prop3.4
\begin{prop}
We set 

\begin{align*}
L&=\partial_t +b\cdot \nabla -\Delta +\frac{2}{r}\partial_r,\\
N&=n\cdot \nabla. 
\end{align*}\\
Assume that $\Gamma\in C^{2,1}(\overline{\Pi}\times (0,T])\cap C(\overline{\Pi}\times [0,T])$ satisfies 

\begin{align*}
(L+1) \Gamma&\leq 0\quad \textrm{in}\ \Pi\times (0,T],   \tag{3.7}\\
(N+2)\Gamma &\leq 0\quad \textrm{on}\ \partial \Pi\times (0,T], \tag{3.8}  \\
\Gamma&\leq 0\quad \textrm{on}\ \Pi\times \{t=0\}.   \tag{3.9} 
\end{align*}\\
Then, 

\begin{align*}
\Gamma\leq 0\quad \textrm{in}\ \Pi\times [0,T].
\end{align*}
\end{prop}

\vspace{5pt}

\begin{proof}[Proof of Proposition 3.2]
Applying Proposition 3.4 for $\tilde{\Gamma}=\Gamma e^{-t}$ implies (3.4).
\end{proof}

\vspace{15pt}

We first consider the case when the function $\Gamma$ attains a maximum in $\overline{\Pi}$. When $\Gamma$ attains the maximum as $|x|\to\infty$, we modify $\Gamma$ so that it attains a maximum in $\overline{\Pi}$.

\vspace{15pt}

\begin{proof}[Proof of Proposition 3.4]
We argue by contradiction. Suppose on the contrary that there exists a point $(x_0,t_0)\in \Pi \times [0,T]$ such that 

\begin{align*}
\Gamma (x_0,t_0)>0.   \tag{3.10}
\end{align*}\\
We set 

\begin{align*}
M=\sup\big\{\Gamma (x,t)\ |\ x\in \Pi,\ t\in [0,T]\ \big\}>0.
\end{align*}\\

\noindent 
\textit{Case 1. The function $\Gamma$ attains the maximum in $\Pi\times [0,T]$.}\\ 

We take a point $(x_1,t_1)\in \Pi\times [0,T]$ such that 

\begin{align*}
M=\Gamma (x_1,t_1)>0.
\end{align*}\\
By (3.9), we may assume that $t_1>0$. Then, there are two cases whether $x_1\in \Pi$ or $x_1\in \partial\Pi$.\\
\noindent 
(a) $x_1\in \Pi$. We observe that 

\begin{align*}
\partial_t \Gamma(x_1,t_1)&\geq 0,\\
\nabla \Gamma (x_1,t_1)&=0,\\
\Delta \Gamma (x_1,t_1)&\leq 0.
\end{align*}\\
Hence we have 

\begin{align*}
((L+1)\Gamma)(x_1,t_1)\geq \Gamma(x_1,t_1)>0. 
\end{align*}\\
This contradicts (3.7). Thus the function $\Gamma$ does not attain the maximum in the interior of $\Pi$.\\
(b) $x_1\in \partial\Pi$. Since the function $\Gamma$ increases along the normal direction near the boundary, we have

\begin{align*}
\frac{\partial \Gamma}{\partial n}(x_1,t_1)\geq 0.
\end{align*} \\
It follows that
 
\begin{align*}
((N+2)\Gamma)(x_1,t_1)\geq 2\Gamma (x_1,t_1)>0.
\end{align*}\\
This contradicts (3.8). Thus the function $\Gamma$ does not attain the maximum on the boundary. \\

\noindent 
\textit{Case 2. The function $\Gamma$ attains the maximum at space infinity.}\\

We modify $\Gamma$ and reduce the problem to Case 1. We set 

\begin{align*}
\Gamma_{\varepsilon}(x,t)=\Gamma (x,t)-\varepsilon (At+|x|^{2}),
\end{align*}\\
by positive constants $A, \varepsilon>0$. We shall show that, by choosing $A^{-1}$ and $\varepsilon$ sufficiently small, depending on $b$, $x_0$, $t_0$ and $\Gamma (x_0,t_0)$, the function $\Gamma_{\varepsilon}$ satisfies the conditions (3.7)-(3.10). Once we verify these conditions, it is not difficult to derive a contradiction. In fact, the function $\Gamma_{\varepsilon}$ is negative in $\overline{\Pi}\cap \{|x|> R\} \times [0,T]$ for $R=\sqrt{M/\varepsilon}$. The condition (3.10) for $\Gamma_{\varepsilon}$ implies the existence of some point $(x_1,t_1)\in \overline{\Pi}\cap \{|x|\leq R \}  \times [0,T]$ such that 

\begin{align*}
M_{\varepsilon}&=\sup\big\{\Gamma_{\varepsilon}(x,t)\ |\ x\in \Pi,\ t\in [0,T]\     \big\} \\
&=\Gamma_{\varepsilon}(x_1,t_1)>0.
\end{align*}\\
However, by the same way as we have shown in Case 1, the conditions (3.7)-(3.10) for $\Gamma_{\varepsilon}$ imply that such the point $(x_1,t_1)$ does not exist. Thus we are able to conclude that Case 2 does not occur neither.

It remains to show (3.7)-(3.10) for $\Gamma_{\varepsilon}$. It follows that 

\begin{align*}
(\partial_n+2)(At+|x|^{2})
&=(-\partial_r +2)(At+r^{2}+|z|^{2})\\
&=2(At+|z|^{2})+2r(r-1)\\
&\geq 0,
\end{align*}
\begin{align*}
(N+2)\Gamma_{\varepsilon}
=(N+2)\Gamma-\varepsilon(\partial_n+2)(At+|x|^{2})
\leq 0.
\end{align*}\\
Thus the conditions (3.8) and (3.9) are satisfied for $A,\varepsilon>0$. We show that (3.7) holds for $\Gamma_{\varepsilon}$ and sufficiently large $A$. Since 

\begin{align*}
L(At+|x|^{2})
&=\Big(\partial_t+ b\cdot \nabla -\Delta +\frac{2}{r}\partial_r\Big)(At+|x|^{2})\\
&=A+2 b\cdot x-2,
\end{align*}\\
it follows that
 
\begin{align*}
(L+1)\Gamma_{\varepsilon}
&=(L+1)\Gamma-\varepsilon (L+1)(At+|x|^{2})\\
&=(L+1)\Gamma-\varepsilon (A(1+t)+|x|^{2}+2 b\cdot x-2). 
\end{align*}\\
Since the function $\Gamma$ satisfies (3.7), the first term of the right-hand side is negative. We set 

\begin{align*}
A_0=\sup\big\{ 2+2||b||_{L^{\infty}(\Pi\times [0,T])}|x|-|x|^{2}\ |\ x\in \Pi\  \big\}>0.
\end{align*} \\
It follows that 

\begin{align*}
A(1+t)+|x|^{2}+2b\cdot x-2
&\geq A-(2+2||b||_{\infty}|x|-|x|^{2})\\
&\geq A-A_0.
\end{align*}\\
Thus the condition (3.7) holds for $\Gamma_{\varepsilon}$ and $A\geq A_0$. Since 

\begin{align*}
\Gamma_{\varepsilon}(x_0,t_0)=\Gamma(x_0,t_0)-\varepsilon(At_0+|x_0|^{2}),
\end{align*}\\
the condition (3.10) holds for $\Gamma_{\varepsilon}$, $\varepsilon<\varepsilon_0$ and $\varepsilon_0=\Gamma(x_0,t_0)(At_0+|x_0|^{2})^{-1}>0$. We proved that (3.7)-(3.10) holds for $\Gamma_{\varepsilon}$. The proof is now complete.
\end{proof}

\vspace{20pt}

\section{An a priori $L^{\infty}$-estimate for swirl}

\vspace{10pt}

We prove the a priori $L^{\infty}$-estimate for the swirl component (1.6) (Lemma 4.7). Since the boundedness of $ru^{\theta}$ does not follow from properties of local-in-time solutions to (1.1), we extend the $L^{\infty}$-estimate (1.9) for mild solutions to (1.8). In the subsequent section, we show that $ru^{\theta}$ is a mild solution to (1.8) and obtain the desired estimate (1.6).

\vspace{15pt}

\subsection{Mild solutions} We define a mild solution of (1.8). We set the elliptic operators by

\begin{align*}
L_0\gamma=\Delta\gamma -\frac{1}{r^{2}}\gamma,\\
L_1 \Gamma =\Delta\Gamma -\frac{2}{r}\partial_r \Gamma,
\end{align*}\\
subject to the Robin boundary conditions, $\partial_n\gamma+\gamma=0$ and  $\partial_n\Gamma+2\Gamma=0$ on $\partial\Pi$. We also set the operator $L_0'=\Delta-r^{-2}$, subject to the Dirichlet boundary condition. By the classical $L^{p}$-estimates for elliptic operators \cite{ADN}, it is known that the operators $B=L_0$, $L_{1}$, $L_0'$ generate $C_0$-analytic semigroups on $L^{p}$ for $p\in (1,\infty)$ \cite[Theorem 3.1.3]{Lunardi}. Moreover, the semigroups are analytic also for $p=\infty$ (see \cite[Corollary 3.1.24]{Lunardi}). By analyticity of the semigroups, they satisfy the regularizing estimate  

\begin{align*}
||\partial_x^{k}e^{tB}f||_{\infty}\leq \frac{C}{t^{\frac{3}{2p}+\frac{|k|}{2}}}||f||_{p}   \tag{4.1}
\end{align*}\\
for $0<t\leq T_0$, $3<p\leq \infty$ and $|k|\leq 1$. By using the semigroup $e^{t L_1}$, we consider the integral equation 

\begin{align*}
\Gamma =e^{tL_1}\Gamma_0-\int_{0}^{t}e^{(t-s)L_1}(b\cdot \nabla \Gamma)(s) \dd s.   \tag{4.2}
\end{align*}\\
We assume that the coefficient $b$ satisfies the regularity condition 

\begin{equation*}
t^{\frac{1}{2}-\frac{3}{2p}}b\in C_{0}([0,T]; L^{p})\quad\textrm{for}\  3<p\leq \infty.  \tag{4.3}
\end{equation*}\\
Here, $C_{0}([0,T]; L^{p})$ denotes the space of all functions in $C([0,T]; L^{p})$, vanishing at time zero. Note that solutions of (1.4) satisfies the condition (4.3) by Lemma 2.1. We prove the $L^{\infty}$-estimate (1.9) for mild solutions $\Gamma\in C_{w}([0,T]; L^{\infty})$ of (4.2), where $C_{w}([0,T]; L^{\infty})$ denotes the space of all weakly-star continuous functions from $[0,T]$ to $L^{\infty}$. \\

We first recall that mild solutions of (4.2) are H\"older continuous up to second orders in $\overline{\Pi} \times [0,T]$ for sufficiently smooth $\Gamma_0$ and $b$ by the H\"older regularity results for second order equations \cite[Chapter IV]{LSU}, \cite[Chapter 5]{Lunardi}. 

Let $C(\overline{\Pi})$ denote the space of all bounded and continuous functions in $\overline{\Pi}$. Let $C^{m}(\overline{\Pi})$ denote the space of all functions $f\in C(\overline{\Pi})$ such that $\partial_x^{k}f\in C(\overline{\Pi})$ for $|k|\leq m$ with non-negative integer $m$. We denote by $C^{\infty}(\overline{\Pi})$ the space of all functions in $C^{m}(\overline{\Pi})$ for all $m\geq 1$. We denote by $C^{\mu}(\overline{\Pi})$ the space of all $\mu$-th H\"older continuous functions $f\in C(\overline{\Pi})$ for $\mu\in (0,1)$. For $m=[m]+\mu$, $C^{m}(\overline{\Pi})$ denotes the space of all functions $f\in C^{[m]}(\overline{\Pi})$ such that $\partial_{x}^{k} f\in C^{\mu}(\overline{\Pi})$ for $|k|=[m]$, where $[m]$ is the greatest integer smaller than $m>0$. We denote by $C^{\mu,\mu/2}(\overline{\Pi}\times [0,T])$ the parabolic H\"older space for $\mu\in (0,2)$, which is the space of all functions $f\in C(\overline{\Pi}\times [0,T])$ such that $f(\cdot,t)\in C^{\mu}(\overline{\Pi})$ for $t\in [0,T]$ and $f(x,\cdot)\in C^{\mu/2}[0,T]$ for $x\in \overline{\Pi}$. We denote by $C^{2+\mu,1+\mu/2}(\overline{\Pi}\times [0,T])$ the space of all functions $f\in C^{2,1}(\overline{\Pi}\times [0,T])$ such that $\partial_t^{s}\partial_x^{k}f\in C^{\mu,\mu/2}(\overline{\Pi}\times [0,T])$ for $2s+|k|\leq 2$.

\vspace{15pt}

%lem4.1
\begin{prop}
Let $T>0$. Let $b$ satisfy (4.3).\\

\noindent
(i) For $\Gamma_0\in L^{\infty}$, there exists a unique mild solution $\Gamma \in C_{w}([0,T]; L^{\infty})$ of (4.2) such that $t^{1/2}\nabla \Gamma \in C_{w}([0,T]; L^{\infty})$. If $\Gamma_0$ and $b$ are axisymmetric, the mild solution $\Gamma$ is axisymmetric. \\

\noindent 
(ii) Assume that 

\begin{align*}
&b\in C^{\mu,\mu/2}(\overline{\Pi}\times [0,T]), \quad \mu\in (0,1),  \tag{4.4}  \\
&\Gamma_0\in C^{2+\mu}(\overline{\Pi})\quad \textrm{and}\quad \partial_n\Gamma+2\Gamma=0\quad \textrm{on}\ \partial\Pi. \tag{4.5}
\end{align*}\\
Then, the mild solution belongs to $C^{2+\mu,1+\mu/2}(\overline{\Pi}\times [0,T])$. In particular, the $L^{\infty}$-estimate (1.9) holds for $t\geq 0$.
\end{prop}

\vspace{5pt}
\begin{proof}
The assertion (i) follows from a standard iteration argument. We are able to prove axial symmetry by a similar way as we did in the proof of Lemma 2.4. The assertion (ii) follows from a H\"older regularity result for second order equations \cite[Theorem 5.1.21, Corollary 5.1.22]{Lunardi}. The $L^{\infty}$-estimate (1.9) follows from Lemma 3.1.  
\end{proof}

\vspace{15pt}

\subsection{Approximation of initial data}
We prove the $L^{\infty}$-estimate (1.9) without the conditions (4.4) and (4.5) by approximation. For this purpose, we prepare H\"older norms for space-time functions \cite{LSU}. We set the $\mu$-th H\"older semi-norm in $Q=\Omega\times (\delta,T]$ for $\mu \in (0,1)$ by

\begin{align*}
&[f]^{(\mu,\frac{\mu}{2})}_{Q}=\sup_{t\in (\delta,T]}[f]^{(\mu)}_{\Omega}(t)+\sup_{x\in \Omega}[f]^{(\frac{\mu}{2})}_{(\delta,T]}(x),\\
&[f]^{(\mu)}_{\Omega}(t)=\sup\Bigg\{\frac{|f(x,t)-f(y,t)|}{|x-y|^{\mu}}\ \Bigg|\ x,y\in \Omega,\ x\neq y  \Bigg\},  \\
&[f]^{(\frac{\mu}{2})}_{(\delta,T]}(x)=\sup\Bigg\{\frac{|f(x,t)-f(x,s)|}{|t-s|^{\frac{\mu}{2}}}\ \Bigg|\ t,s\in (\delta,T],\ t\neq s  \Bigg\}.
\end{align*}\\
When $\mu=1$, we set

\begin{align*}
[f]^{(1,\frac{1}{2})}_{Q}=||\nabla f||_{L^{\infty}(Q)}+\sup_{x\in \Omega}[f]_{(\delta,T]}^{(\frac{1}{2})}(x).
\end{align*}\\
For $m=[m]+\mu$, we set  
 
\begin{align*}
&[f]_{Q}^{(m,\frac{m}{2})}=\sum_{2s+|k|=[m]}[\partial_t^{s}\partial_x^{k}f]^{(\mu,\frac{\mu}{2})}_{Q},\\
&|f|_{Q}^{(m,\frac{m}{2})}=\sum_{2s+|k|\leq [m]}||\partial_t^{s}\partial_x^{k}f||_{L^{\infty}(Q)}+[f]_{Q}^{(m,\frac{m}{2})}.
\end{align*}

\vspace{15pt}

We first remove the condition (4.5) by approximation of $\Gamma_0\in L^{\infty}$.

\vspace{15pt}

%prop4.2
\begin{prop}
For $\Gamma_0\in L^{\infty}(\Pi)$, there exists a sequence $\{\Gamma_{0,\varepsilon}\}\subset  C^{\infty}(\overline{\Pi})$ supported in $\Pi$ such that 

\begin{equation*}
\begin{aligned}
&||\Gamma_{0,\varepsilon}||_{\infty}\leq ||\Gamma_0||_{\infty}\\
&\Gamma_{0,\varepsilon}\to \Gamma_0\quad \textrm{a.e.}\ \textrm{in}\ \Pi.
\end{aligned}
\tag{4.6}
\end{equation*}
\end{prop}

\vspace{5pt}

\begin{proof}
For $x=re_{r}(\theta)+ze_{z}$, we set 

\begin{equation*}
\tilde{\Gamma}_{0,\varepsilon}(x)=
\begin{cases}
&\Gamma_0((r-\varepsilon)e_{r}(\theta)+ze_{z})\hspace{13pt} r\geq 1+\varepsilon,\\
&0\hspace{100pt} 0\leq r<1+\varepsilon.
\end{cases}
\end{equation*}\\
By mollification of $\tilde{\Gamma}_{0, \varepsilon}$, we obtain the desired sequence.
\end{proof}

\vspace{15pt}

%prop4.3
\begin{prop}
In Proposition 4.1 (ii), the estimate (1.9) holds without the condition (4.5).
\end{prop}

\vspace{5pt}

\begin{proof}
For $\Gamma_0\in L^{\infty}$, we take a sequence $\{\Gamma_{0,\varepsilon}\}$ satisfying (4.6). Since $\Gamma_{0,\varepsilon}$ is smooth in $\overline{\Pi}$ and supported in $\Pi$, it satisfies the condition (4.5). Since the estimate (1.9) holds for the mild solution $\Gamma_{\varepsilon}$ of (4.2) for $\Gamma_{0,\varepsilon}$ by Proposition 4.1 (ii), it follows from (4.6) that 

\begin{align*}
||\Gamma_{\varepsilon}||_{\infty}\leq ||\Gamma_{0}||_{\infty}\quad t>0.  \tag{4.7}
\end{align*}\\
We shall show that $\Gamma_{\varepsilon}$ converges to a mild solution of (4.2) for $\Gamma_0$. We use the H\"older continuity of the coefficient $b$ in (4.4). We apply the local H\"older estimate for parabolic equations \cite[Chapter IV, Theorem 10.1]{LSU} and estimate

\begin{align*}
|\Gamma_{\varepsilon}|^{(2+\mu,1+\frac{\mu}{2})}_{Q}\leq C||\Gamma_{\varepsilon}||_{L^{\infty}(\Pi\times (0,T))}    \tag{4.8}
\end{align*}\\
for $Q=(B\cap \Pi)\times (\delta,T]$ and $\delta >0$ with some constant $C$, independent of $\varepsilon$. Here, $B\subset \mathbb{R}^{3}$ denotes an open ball satisfying $B\cap \Pi\neq \emptyset$. By (4.7) and (4.8), $\Gamma_{\varepsilon}$ subsequently converges to a limit $\Gamma$ locally uniformly in $\overline{\Pi}\times (0,T]$ up to second derivatives. 

It is not difficult to see that the limit $\Gamma$ is a mild solution of (4.2) for $\Gamma_0$. In fact, by choosing a subsequence, we have 

\begin{align*}
e^{tL_1}\Gamma_{0,\varepsilon}\to e^{tL_1}\Gamma_0\quad \textrm{locally uniformly in}\ \overline{\Pi}\times (0,T].
\end{align*}\\
Since $\nabla \Gamma_{\varepsilon}$ converges to $\nabla \Gamma$ locally uniformly in $\overline{\Pi}\times (0,T]$, similarly for each $0<s<t$, we have 

\begin{align*}
e^{\rho L_1}b\cdot \nabla \Gamma_{\varepsilon}\to e^{\rho L_1}b\cdot \nabla \Gamma\quad \textrm{locally uniformly in}\ \overline{\Pi}\times (0,T].
\end{align*}\\
Hence sending $\varepsilon\to 0$ to (4.2) implies the limit $\Gamma$ is a mild solution for $\Gamma_0$. The estimate (4.7) is inherited to the limit $\Gamma$.
\end{proof}

\vspace{15pt}

\subsection{Approximation of a coefficient}

\vspace{15pt}
We next remove the condition (4.4).

\vspace{15pt}

%prop4.4
\begin{prop}
For $b$ satisfying (4.3), there exists a sequence $\{b_{\varepsilon}\}\subset C^{\infty}(\overline{\Pi}\times [0,T])$ satisfying (4.3) and

\begin{align*}
\lim_{\varepsilon\to0}\sup_{0\leq t\leq T}t^{\frac{1}{2}-\frac{3}{2p}}||b-b_{\varepsilon}||_{p}(t)=0\quad \textrm{for}\ 3<p< \infty. \tag{4.9}
\end{align*}
\end{prop}

\vspace{5pt}

\begin{proof}
We may assume that $b$ is smooth in $\overline{\Pi}$ by mollification by spatial variables. Since $g=t^{1/2-3/(2p)}b$ vanishes at time zero by (4.3), by shifting $g$ by a time variable, and mollification, we obtain a sequence $\{g_{\varepsilon}\}\subset C^{\infty}(\overline{\Pi}\times [0,T])$ such that $g_{\varepsilon}(\cdot ,t)$ is supported in $(0,T]$ and 

\begin{align*}
\lim_{\varepsilon\to 0}\sup_{0\leq t\leq T}||g_{\varepsilon}-g||_{p}(t)=0.   \tag{4.10}
\end{align*}\\
Since $g_{\varepsilon}(\cdot, t)$ is supported in $(0,T]$, the function $b_{\varepsilon}=t^{-1/2+3/2p}g_{\varepsilon}$ is smooth in $\overline{\Pi}\times [0,T]$ and satisfies (4.3). The convergence (4.9) follows from (4.10).
\end{proof}

\vspace{5pt}

%lem4.5
\begin{lem}
The estimate (1.9) holds for mild solutions of (4.2) for $\Gamma_0\in L^{\infty}$ and $t>0$.
\end{lem}

\vspace{5pt}

\begin{proof}
We shall show the estimate (1.9) between $0<t\leq T_{1}$ for some $T_1>0$. Once we have (1.9) near time zero, it is extendable for all $t>0$ by taking $t=T_1$ as an initial time. We take a sequence $\{b_{\varepsilon}\}$ satisfying (4.9). Since the estimate (1.9) holds for a mild solution $\Gamma_{\varepsilon}$ for $\Gamma_0\in L^{\infty}$ and the coefficient $b_{\varepsilon}$ by Proposition 4.3, we have 

\begin{align*}
||\Gamma_{\varepsilon}||_{\infty}\leq ||\Gamma_0||_{\infty}\quad t>0.  \tag{4.11}
\end{align*}\\
We shall show that $\Gamma_{\varepsilon}$ converges to a mild solution $\Gamma$ in the sense that 

\begin{align*}
\lim_{\varepsilon\to0}\sup_{0< t\leq T_1}\Big\{||\Gamma-\Gamma_{\varepsilon}||_{\infty}
+t^{\frac{1}{2}}||\nabla (\Gamma-\Gamma_{\varepsilon})||_{\infty}
\Big\}=0. \tag{4.12}
\end{align*}\\
The desired estimate follows from (4.11) and (4.12) by sending $\varepsilon\to 0$.\\

We set $\rho_{\varepsilon}=\Gamma-\Gamma_{\varepsilon}$ and $a_{\varepsilon}=b-b_{\varepsilon}$. It follows from (4.2) that  

\begin{align*}
\rho_{\varepsilon}
&=-\int_{0}^{t}e^{(t-s)L_1}(b\cdot \nabla \Gamma-b_{\varepsilon}\cdot \nabla \Gamma_{\varepsilon})\dd s\\
&=-\int_{0}^{t}e^{(t-s)L_1}(a_{\varepsilon}\cdot \nabla \Gamma+b_{\varepsilon}\cdot \nabla \rho_{\varepsilon})\dd s.
\end{align*}\\
For $p\in (3,\infty)$, we set the constants

\begin{align*}
K_{\varepsilon}&=\sup_{0< t\leq T_1}\big\{||\rho_{\varepsilon}||_{\infty}+t^{\frac{1}{2}}||\nabla \rho_{\varepsilon}||_{\infty}   \big\},\\
K&=\sup_{0< t\leq T_1}\big\{||\Gamma||_{\infty}+t^{\frac{1}{2}}||\nabla \Gamma||_{\infty}   \big\},\\
N_{\varepsilon}&=\sup_{0\leq t\leq T_1}t^{\frac{1}{2}-\frac{3}{2p} }||a_{\varepsilon}||_{p},\\
L_{\varepsilon}&=\sup_{0\leq t\leq T_1}t^{\frac{1}{2}-\frac{3}{2p} }||b_{\varepsilon}||_{p}.
\end{align*}\\
It follows from (4.1) that

\begin{align*}
||\rho_{\varepsilon}||_{\infty}
&\leq \int_{0}^{t}\frac{C}{(t-s)^{\frac{3}{2p}}}\big(||a_{\varepsilon}||_{p}||\nabla \Gamma||_{\infty}+||b_{\varepsilon}||_{p}||\nabla \rho_{\varepsilon}||_{\infty}  \big)\dd s\\
&\leq C(N_{\varepsilon}K+L_{\varepsilon}K_{\varepsilon})\int_{0}^{t}\frac{\dd s}{(t-s)^{\frac{3}{2p}}s^{1-\frac{3}{2p}}  }\\
&=C_0(N_{\varepsilon}K+L_{\varepsilon}K_{\varepsilon}).
\end{align*}\\
Similarly, we estimate $\nabla \rho_{\varepsilon}$ and obtain 

\begin{align*}
K_{\varepsilon}\leq C_0(N_{\varepsilon}K+L_{\varepsilon}K_{\varepsilon}).
\end{align*}\\
We take an arbitrary $\delta>0$. By (4.3), there exists $T_1>0$ such that 

\begin{align*}
\sup_{0<t\leq T_1}t^{\frac{1}{2}-\frac{3}{2p}}||b||_{p}\leq \delta.
\end{align*}\\
By (4.9), that there exits $\varepsilon_0>0$ such that 

\begin{align*}
\sup_{0<t\leq T_1}t^{\frac{1}{2}-\frac{3}{2p}}||b-b_{\varepsilon}||_{p}\leq \delta\quad \textrm{for}\ \varepsilon\leq \varepsilon_0.
\end{align*}\\
We estimate 

\begin{align*}
L_{\varepsilon}=\sup_{0<t\leq T_1}t^{\frac{1}{2}-\frac{3}{2p}}||b_{\varepsilon}||_{p}
&\leq \sup_{0<t\leq T_1}t^{\frac{1}{2}-\frac{3}{2p}}||b_{\varepsilon}-b||_{p}+\sup_{0<t\leq T_1}t^{\frac{1}{2}-\frac{3}{2p}}||b||_{p}\\
&\leq 2\delta\qquad \textrm{for}\ \varepsilon\leq \varepsilon_0.
\end{align*}\\
By taking $\delta=(4C_0)^{-1}$, we estimate 

\begin{align*}
K_{\varepsilon}\leq 2C_0N_{\varepsilon}K.
\end{align*}\\
Since $N_{\varepsilon}\to 0$ as $\varepsilon \to0$ by (4.9), we proved (4.12).
\end{proof}

\vspace{15pt}

%subsection 4.3
\subsection{An application to axisymmetric solutions}

We now prove the a priori $L^{\infty}$-estimate for the swirl component (1.6). It suffices to show that the swirl component $ru^{\theta}$ is a mild solution of (4.2).

\vspace{15pt}

%prop4.9
\begin{prop}
The semigroups satisfy 
 
\begin{align*}
e_{\theta}\cdot e^{t A}f&=e^{t L_0}f^{\theta},  \tag{4.13}\\
e_{\theta}\cdot \textrm{curl}\ e^{tA}g&=e^{tL_{0}'}\big(e_{\theta}\cdot \textrm{curl}\ g \big),    \tag{4.14}   \\
r e^{tL_0}\gamma &=e^{tL_1}(r\gamma),  \tag{4.15}
\end{align*}\\
for axisymmetric $f$, $g\in L^{2}_{\sigma}$ and $\gamma\in L^{\infty}$ satisfying $e_{\theta}\cdot \textrm{curl}\ g\in L^{2}$ and $r\gamma \in L^{\infty}$.
\end{prop}

\vspace{5pt}

\begin{proof}
We set $w=e^{tA}f$. Since $w^{\theta}=e_{\theta}\cdot w$ satisfies 

\begin{align*}
\partial_t w^{\theta}-\Big(\Delta-\frac{1}{r^{2}}\Big)w^{\theta}&=0\quad \textrm{in}\ \Pi\times (0,\infty),\\
\partial_n w^{\theta}+w^{\theta}&=0\quad \textrm{on}\ \partial\Pi\times (0,\infty),\\
w^{\theta}&=f^{\theta}\hspace{6pt} \textrm{on}\ \Pi\times \{t=0\},
\end{align*}\\
the function $w^{\theta}$ agrees with $e^{tL_0}f^{\theta}$ by the uniqueness of the heat equation. Similarly, we are able to prove (4.14) and (4.15).
\end{proof}

\vspace{15pt}

%lem 4.8
\begin{lem}[A priori $L^{\infty}$-estimate]
Let $u=v+u^{\theta}e_{\theta}$ be an axisymmetric mild solution of (1.4) in Lemma 2.4. Assume that $ru^{\theta}_{0}\in L^{\infty}$. Then, $ru^{\theta}$ is a mild solution of (4.2) for $b=v$. In particular, the $L^{\infty}$-estimate (1.6) holds. 
\end{lem}

\vspace{5pt}

\begin{proof}
We set 

\begin{equation*}
\begin{aligned}
h
&=u\cdot \nabla u
=\Big(v\cdot \nabla u^{r}-\frac{|u^{\theta}|^{2}}{r}\Big)e_{r}
+\Big(v\cdot \nabla u^{\theta}+\frac{u^{r}}{r}u^{\theta} \Big)e_{\theta}
+(v\cdot \nabla u^{z})e_{z},\\
\nabla \Phi&=(I-\mathbb{P})h,
\end{aligned}
\tag{4.16}
\end{equation*}\\
for axisymmetric mild solutions $u$ in $[0,T]$. Since $h$ is axisymmetric, the function $\Phi$ is independent of $\theta$ and we have 

\begin{align*}
e_{\theta}\cdot \mathbb{P}u\cdot \nabla u
&=e_{\theta}\cdot(h-\nabla \Phi) \\
&=e_{\theta}\cdot(h-e_{r}\partial_r \Phi-e_{z}\partial_z \Phi) \\
&=v\cdot \nabla u^{\theta}+\frac{u^{r}}{r}u^{\theta}.
\end{align*}\\
We multiply $e_{\theta}$ by (1.4). It follows from (4.13) that 

\begin{align*}
u^{\theta}=e^{tL_0}u^{\theta}_0-\int_{0}^{t}e^{(t-s)L_0}\Big(v\cdot \nabla u^{\theta}+\frac{u^{r}}{r}u^{\theta}\Big)\dd s.
\end{align*}\\
Since $u^{\theta}_{0}\in L^{\infty}$ by $ru^{\theta}_{0}\in L^{\infty}$, the above integral form implies that $u^{\theta}\in C_{w}([0,T]; L^{\infty})$ and $t^{1/2}\nabla u^{\theta}\in C_{w}([0,T]; L^{\infty})$.

On the other hand, there exists a unique axisymmetric mild solution $\Gamma$ for $\Gamma_0=ru_{0}^{\theta}$ and $b=v$ by Proposition 4.1 (i). We multiply $r^{-1}$ by (4.2). It follows from (4.15) that $\gamma=\Gamma/r$ satisfies 

\begin{align*}
\gamma=e^{tL_0}u_0^{\theta}-\int_{0}^{t}e^{(t-s)L_0}\Big(v\cdot \nabla \gamma+\frac{u^{r}}{r}\gamma\Big)\dd s.
\end{align*}\\
Since $\gamma\in C_{w}([0,T]; L^{\infty})$ and $t^{1/2}\nabla \gamma \in C_{w}([0,T]; L^{\infty})$, it is not difficult to show that $u^{\theta}$ agrees with $\gamma$ by estimating the difference $u^{\theta}-\gamma$. Thus $r u^{\theta}$ is a mild solution of (4.2). The proof is now complete.
\end{proof}

\vspace{15pt}

%セクション5
\section{Energy estimates for the azimuthal component of vorticity}
\vspace{10pt}

We prove the global estimates (1.11) and (1.12). 

\begin{equation*}
\begin{aligned}
\partial_t \omega^{\theta} +v\cdot \nabla \omega^{\theta}-\frac{u^{r}}{r}\omega^{\theta} -\Big(\Delta -\frac{1}{r^{2}}\Big)\omega^{\theta} &=\frac{\partial_z |u^{\theta}|^{2}}{r}\quad \textrm{in}\ \Pi\times (0,T), \\
\omega^{\theta}&=0\hspace{35pt} \textrm{on}\ \partial\Pi\times (0,T),\\
\omega^{\theta}&=\omega^{\theta}_0\hspace{29pt} \textrm{on}\ \Pi\times \{t=0\}.
\end{aligned}
\tag{5.1}
\end{equation*}\\

%prop5.1
\begin{prop}
Axisymmetric mild solutions of (1.4) in Lemma 2.4 satisfy 

\begin{align*}
\omega^{\theta}\in C^{\gamma}((0,T]; D(L_0') )\cap C^{1+\gamma}((0,T]; L^{2} ),\quad 0<\gamma<\frac{1}{2}. \tag{5.2}
\end{align*}\\
In particular, $\omega^{\theta}$ satisfies the vorticity equation (5.1), where $D(L_0')=H^{2}\cap H^{1}_{0}$ denotes the domain of the operator $L_0'=\Delta-r^{-2}$ on $L^{2}$.
\end{prop}

\vspace{5pt}

\begin{proof}
We recall that the mild solution $u$ is expressed by

\begin{align*}
u(\eta+\delta)=e^{\eta A}u(\delta)+\int_{0}^{\eta}e^{(\eta-s)A}f(\delta+s)\dd s,
\end{align*}\\
for $0\leq \delta \leq \eta\leq T-\delta$ and $f=-\mathbb{P}u\cdot \nabla u$. It follows from (4.16) that 

\begin{align*}
g:&=e_{\theta}\cdot \textrm{curl}\ f\\
&=-e_{\theta}\cdot \textrm{curl}\ (h-\nabla \Phi)\\
&=-\partial_zh^{r}+\partial_rh^{z}\\
&=-v\cdot \nabla \omega^{\theta}+\frac{u^{r}}{r}\omega^{\theta}+\frac{\partial_z|u^{\theta}|^{2}}{r}.
\end{align*}\\
We multiply $e_{\theta}\cdot \textrm{curl}$ by $u$. It follows from (4.14) that 

\begin{align*}
\omega^{\theta}(\eta+\delta)=e^{\eta L_0'}\omega^{\theta}(\delta)+\int_{0}^{\eta}e^{(\eta-s)L_0'}g(\delta+s)\dd s.
\end{align*}\\
Since $u\in C^{\gamma}((0,T]; H^{2} )$ for $\gamma\in (0,1/2)$ by (2.4), we have $g\in C^{\gamma}((0,T]; L^{2})$. By Proposition 2.2, $\omega^{\theta}$ satisfies (5.2).
\end{proof}

\vspace{15pt}

%prop5.2
\begin{prop}
The function $v=u^{r}e_{r}+u^{z}e_{z}$ satisfies 

\begin{align*}
||\nabla v||_{2}&=||\omega^{\theta}||_{2},  \tag{5.3} \\
||u^{r}||_{4}&\leq C ||u^{r}||_{2}^{\frac{1}{4}}||\omega^{\theta}||_{2}^{\frac{3}{4}},   \tag{5.4} \\
||u^{z}||_{4}&\leq C ||u^{z}||_{2}^{\frac{1}{4}}(||u^{z}||_{2}+||\omega^{\theta}||_{2})^{\frac{3}{4}},   \tag{5.5} \\
||\omega^{\theta}||_{4}&\leq C ||\omega^{\theta}||_{2}^{\frac{1}{4}}||\nabla \omega^{\theta}||_{2}^{\frac{3}{4}}, \tag{5.6} \\
||\nabla (\omega^{\theta}e_{\theta})||_{2}&=\Big(||\nabla\omega^{\theta}||_{2}^{2}+\Big\|\frac{\omega^{\theta}}{r} \Big\|_{2}^{2}\Big)^{\frac{1}{2}},  \tag{5.7}
\end{align*}\\
with some constant $C$. In particular, the estimate (1.10) holds.
\end{prop}

\vspace{5pt}

\begin{proof}
Since $v$ satisfies 
 
\begin{align*}
-\Delta v
&=\textrm{curl}\ \textrm{curl}\ v-\nabla \textrm{div}\ v \\
&=\textrm{curl}\ (\omega^{\theta}e_{\theta}),
\end{align*}\\
it follows from (2.6) that 

\begin{align*}
\int_{\Pi}|\nabla v|^{2}\dd x
&=-\int_{\Pi}\Delta v \cdot v\dd x
+\int_{\partial\Pi}\frac{\partial v}{\partial n}\cdot v\dd{ \mathcal{H}}\\
&=\int_{\Pi}\textrm{curl}\ (\omega^{\theta}e_{\theta})\cdot v\dd x
-\int_{\partial\Pi}(\partial_r u^{r} u^{r}+\partial_{r} u^{z} u^{z})\dd{ \mathcal{H}}\\
&=\int_{\Pi}|\omega^{\theta}|^{2}\dd x.
\end{align*}\\
Thus (5.3) holds. Similarly, we obtain (5.7) by integration by parts. Since $v^{r}$ and $\omega^{\theta}$ vanish on $\partial\Pi$, applying the interpolation inequality (B.1) implies (5.4) and (5.6). We apply (B.2) for $v^{z}$ and obtain (5.5).
\end{proof}

\vspace{15pt}

%lem5.3
\begin{lem}
The estimates (1.11) and (1.12) hold for $t>0$ for axisymmetric mild solutions for $u_0\in \tilde{L}^{3}_{\sigma}$ satisfying $ru_0^{\theta}\in L^{\infty}$ and $\nabla u_0\in L^{2}$ with some constant C. 
\end{lem}

\vspace{5pt}

\begin{proof}
We prove (1.11). Since $\omega^{\theta}/r$ satisfies (1.2) and vanishes on $\partial\Pi$, by multiplying $2\omega^{\theta}/r$ by (1.2) and integration by parts, we have

\begin{align*}
\frac{\dd}{\dd t} \int_{\Pi}\Big|\frac{\omega^{\theta}}{r}\Big|^{2}\dd x
+2\int_{\Pi}\Big|\nabla \Big(\frac{\omega^{\theta}}{r}\Big)\Big|^{2}\dd x =-2\int_{\Pi}\Big(\frac{u^{\theta}}{r}\Big)^{2}\partial_z \Big(\frac{\omega^{\theta}}{r}\Big)\dd x.
\end{align*}\\
Since $r\geq 1$, it follows that 

\begin{align*}
\Big\|\frac{u^{\theta}}{r}\Big\|_{4}
\leq ||u^{\theta}||_{4} 
=\Big\|(ru^{\theta})^{\frac{1}{2}}\Big(\frac{u^{\theta}}{r}\Big)^{\frac{1}{2}}\Big\|_{4}
\leq ||ru^{\theta}||_{\infty}^{\frac{1}{2}}
\Big\|\frac{u^{\theta}}{r}\Big\|_{2}^{\frac{1}{2}}.   \tag{5.8}
\end{align*}\\
By the Young's inequality, we estimate 

\begin{align*}
\Big|2\int_{\Pi}\Big(\frac{u^{\theta}}{r}\Big)^{2}\partial_z \Big(\frac{\omega^{\theta}}{r}\Big)\dd x\Big|
&\leq \Big\|\frac{u^{\theta}}{r}\Big\|_{4}^{4}
+\Big\|\nabla \Big(\frac{\omega^{\theta}}{r}\Big)\Big\|_{2}^{2} \\
&\leq ||ru^{\theta}||_{\infty}^{2}
\Big\|\frac{u^{\theta}}{r}\Big\|_{2}^{2}
+\Big\|\nabla \Big(\frac{\omega^{\theta}}{r}\Big)\Big\|_{2}^{2}.
\end{align*}\\
Hence

\begin{align*}
\frac{\dd}{\dd t} \int_{\Pi}\Big|\frac{\omega^{\theta}}{r}\Big|^{2}\dd x
+\int_{\Pi}\Big|\nabla \Big(\frac{\omega^{\theta}}{r}\Big)  \Big|^{2}\dd x
\leq ||ru^{\theta}||_{\infty}^{2}\Big\|\frac{u^{\theta}}{r}\Big\|_{2}^{2}.
\end{align*}\\
We integrate the both sides between $(0,t)$. By applying (1.5) and (1.6), we obtain (1.11).\\

We prove (1.12). We multiply $2\omega^{\theta}$ by (5.1) to see that 

\begin{align*}
\frac{\dd}{\dd t} \int_{\Pi}|\omega^{\theta}|^{2}\dd x
+2\int_{\Pi}\Big(|\nabla \omega^{\theta}|^{2}+\Big|\frac{\omega^{\theta}}{r}\Big|^{2}\Big)\dd x 
&=2\int_{\Pi}\frac{u^{r}}{r}\omega^{\theta}\omega^{\theta} \dd x
+2\int_{\Pi}\frac{\partial_z|u^{\theta}|^{2}}{r} \omega^{\theta}\dd x\\
&=: I+II.
\end{align*}\\
It follows from (5.4), (5.6) and (5.7) that 

\begin{align*}
|I|=\Big|2\int_{\Pi}\frac{u^{r}}{r}\omega^{\theta}\omega^{\theta}\dd x  \Big|
&\leq 2\Big\|\frac{\omega^{\theta}}{r}\Big\|_{2}||u^{r}||_{4}||\omega^{\theta}||_{4}\\
&\leq C\Big\|\frac{\omega^{\theta}}{r}\Big\|_{2}||u^{r}||_{2}^{\frac{1}{4}}||\omega^{\theta}||_{2}||\nabla \omega^{\theta}||_{2}^{\frac{3}{4}}\\
&\leq C\Big\|\frac{\omega^{\theta}}{r}\Big\|_{2}^{\frac{3}{4}}||u^{r}||_{2}^{\frac{1}{4}}||\omega^{\theta}||_{2}||\nabla( \omega^{\theta}e_{\theta})||_{2}\\
&\leq C'\Big\|\frac{\omega^{\theta}}{r}\Big\|_{2}^{\frac{3}{2}}||u^{r}||_{2}^{\frac{1}{2}}||\omega^{\theta}||_{2}^{2}+\frac{1}{2}||\nabla( \omega^{\theta}e_{\theta})||_{2}^{2}.
\end{align*}\\
Since $r\geq 1$, it follows from (5.8) that 

\begin{align*}
|II|=\Big|2\int_{\Pi}\frac{\partial_z|u^{\theta}|^{2}}{r}\omega^{\theta}\dd x\Big|
&\leq 2||u^{\theta}||_{4}^{2}||\nabla \omega^{\theta}||_{2}\\
&\leq 2||ru^{\theta}||_{\infty}\Big\|\frac{u^{\theta}}{r}\Big\|_{2}||\nabla (\omega^{\theta}e_{\theta})||_{2}\\
&\leq 2||ru^{\theta}||_{\infty}^{2}\Big\|\frac{u^{\theta}}{r}\Big\|_{2}^{2}
+\frac{1}{2}||\nabla (\omega^{\theta}e_{\theta})||_{2}^{2}.
\end{align*}\\
By combining the estimates for $I$ and $II$, we obtain 

\begin{align*}
\frac{\dd}{\dd t} \int_{\Pi}|\omega^{\theta}|^{2}\dd x
+\int_{\Pi}\Big(|\nabla \omega^{\theta}|^{2}+\Big|\frac{\omega^{\theta}}{r}\Big|^{2}\Big)\dd x 
\leq C\Big(\Big\|\frac{\omega^{\theta}}{r}\Big\|_{2}^{\frac{3}{2}}||u^{r}||_{2}^{\frac{1}{2}}+||ru^{\theta}||_{\infty}^{2}  \Big)\Big(||\omega^{\theta}||_{2}^{2}+\Big\|\frac{u^{\theta}}{r}\Big\|_{2}^{2}\Big).
\end{align*}\\
We integrate the both sides between $(0,t)$. By (1.11), (1.5) and (1.6), we obtain (1.12). 
\end{proof}

\vspace{20pt}

%セクション6
\section{Global bounds on $L^{4}$}
\vspace{10pt}

%定理1.1の証明
\begin{proof}[Proof of Theorem 1.1]
For an axisymmetric $u_0\in \tilde{L}^{3}_{\sigma}$ satisfying $ru^{\theta}_{0}\in L^{\infty}$, the axisymmetric mild solution $u\in C([0,T]; \tilde{L}^{3})$ satisfies (1.6) by Lemma 4.7. It follows from (1.6) and (2.2) that $ru^{\theta}_{0}(\cdot , t_0)\in L^{\infty}$ and 

\begin{align*}
\nabla u(\cdot , t_0)\in L^{2}\quad \textrm{for}\ t_0\in (0,T].
\end{align*}\\
We may assume that $\nabla u_0\in L^{2}$ by taking $t=t_0$ as an initial time. It follows from (1.7) and (1.10)-(1.12) that $u\in L^{\infty}(0,\infty; L^{4})$. By (1.5) and Lemma 2.1, the mild solution belongs to $BC([0,\infty); \tilde{L}^{3})$. The proof is now complete.
\end{proof}

\vspace{15pt}

\begin{rems}
(i) (The Euler equations) We constructed global solutions in the exterior domain by using viscosity. For the Euler equations, existence of global solutions is unknown. We refer to \cite{HL13a}, \cite{HL13b} for a one-dimensional blow-up model of axisymmetric Euler flows on the boundary. See \cite{CKY}, \cite{CHKLSY} for blow-up results of models.

\noindent
(ii) (The Dirichlet boundary condition)
The statement of Lemma 2.1 is valid also for the Dirichlet boundary condition. However, in this case unique existence of global solution is unknown even for axisymmetric data without swirl. Since the azimuthal component of vorticity $\omega^{\theta}$ does not vanish on the boundary subject to the Dirichlet boundary condition, the global vorticity estimates (1.11) and (1.12) are not available unlike the slip boundary condition.

\noindent
(iii) (Uniform estimates)
The assertion of Theorem 1.1 is valid also for the exterior of the cylinder $\Pi^{\varepsilon}=\{r>\varepsilon\}$ and we are able to construct global solutions $u=u_{\varepsilon}$ satisfying (1.6) and the energy equality 

\begin{equation*}
\begin{aligned}
\int_{\Pi^{\varepsilon}}|u|^{2}\dd x
&+2\int_{0}^{t}\hspace{-3pt}\int_{\Pi^{\varepsilon}}\Big(|\nabla v|^{2}+|\nabla u^{\theta}|^{2}
+\Big|\frac{u^{\theta}}{r}\Big|^{2}\Big)\dd x\dd s\\
&+\frac{2}{\varepsilon} \int_{0}^{t}\hspace{-3pt}\int_{\partial \Pi^{\varepsilon}}|u^{\theta}|^{2}\dd {\mathcal{H}} \dd s= \int_{\Pi^{\varepsilon}}|u_{0}|^{2}\dd x,   \quad t\geq 0.
\end{aligned}
\tag{6.1}
\end{equation*}\\
For the case without swirl, the a priori estimates

\begin{align*}
\int_{\Pi^{\varepsilon}}\Big|\frac{\omega^{\theta}}{r}\Big|^{2}\dd x
+2\int_{0}^{t}\int_{\Pi^{\varepsilon}}\Big|\nabla \Big(\frac{\omega^{\theta}}{r}\Big)\Big|^{2}\dd x\dd s
\leq \int_{\Pi^{\varepsilon}}\Big|\frac{\omega^{\theta}_{0}}{r}\Big|^{2}\dd x,  \tag{6.2}
\end{align*}

\begin{equation*}
\begin{aligned}
\int_{\Pi^{\varepsilon}}|\omega^{\theta}|^{2}\dd x
&+\int_{0}^{t}\int_{\Pi^{\varepsilon}}\Big(|\nabla\omega^{\theta}|^{2}+\Big|\frac{\omega^{\theta}}{r}\Big|^{2}   \Big)\dd x\dd s \\
&\leq \int_{\Pi^{\varepsilon}}|\omega^{\theta}_{0}|^{2}\dd x  
+C\Big\|\frac{\omega_{0}^{\theta}}{r}\Big\|_{L^{2}(\Pi^{\varepsilon})}^{\frac{3}{2}} ||u_0||_{L^{2}(\Pi^{\varepsilon})}^{\frac{5}{2}},  \quad t>0,
\end{aligned}
\tag{6.3}
\end{equation*}\\
hold with some constant $C$, independent of $\varepsilon$. Hence we have a uniform bound

\begin{align*}
\sup_{\varepsilon\leq \varepsilon_0}||u||_{L^{\infty}(0,\infty; H^{1}(\Pi^{\varepsilon}) )}<\infty,
\end{align*}\\
provided that $L^{2}$-norms of $u^{\theta}_{0}$, $\omega^{\theta}_{0}/r$ and $\omega^{\theta}_{0}$ in $\Pi^{\varepsilon}$ are uniformly bounded for $\varepsilon\leq \varepsilon_0$.
\end{rems}

\vspace{20pt}

\appendix

%Appendix A
\section{Mild solutions on $\tilde{L}^{3}$}

\vspace{10pt}

We give a proof for local solvability of (1.1) on $\tilde{L}^{3}$ (Lemma 2.1). 

\vspace{10pt}

%prop A.1
\begin{prop}
The Stokes semigroup satisfies 

\begin{align*}
||\partial_x^{k}e^{tA}f||_{\tilde{L}^{p}}&\leq \frac{C}{t^{\frac{3}{2}(\frac{1}{q}-\frac{1}{p})+\frac{|k|}{2}}}||f||_{\tilde{L}^{q}},     \tag{A.1}\\
||\partial_x^{k}(e^{\rho A}-1)e^{\eta A}f||_{\tilde{L}^{q}}
&\leq C\frac{\rho^{\alpha}}{\eta^{\alpha+\frac{|k|}{2}}}||f||_{\tilde{L}^{q}},\tag{A.2}
\end{align*}\\
for $f\in \tilde{L}^{q}_{\sigma}$, $2\leq q\leq p< \infty$, $0<t,\rho,\eta\leq T_0$, $\alpha\in (0,1)$, $|k|\leq 1$ and $T_0>0$.
\end{prop}

\vspace{5pt}

\begin{proof}
The estimates (A.1) and (A.2) follow from estimates of the Stokes semigroup on $\tilde{L}^{q}$ \cite[Theorem 1.2]{FR} and the interpolation inequality (B.2).  
\end{proof}

\vspace{10pt}

%propA.2
\begin{prop}
For $u_0\in \tilde{L}^{3}_{\sigma}$, there exists $T>0$ and a unique mild solution of (1.4) satisfying (2.1) and (2.2).
\end{prop}

\vspace{5pt}

\begin{proof}
We set 
\begin{align*}
u_{j+1}&=e^{tA}u_0-\int_{0}^{t}e^{(t-s)A}\mathbb{P}(u_{j}\cdot \nabla u_j)\dd s,\\
u_1&=e^{tA}u_0,
\end{align*}
\begin{align*}
K_{j}=\sup_{0\leq t\leq T}t^{\gamma}\{||u_j||_{\tilde{L}^{p}}+t^{\frac{1}{2}}||\nabla u_j||_{\tilde{L}^{p}} \},
\end{align*}\\
for $\gamma=1/2-3/(2p)$ and $p\in (3,\infty)$. We take $q\in [2,p]$. By applying the Young's inequality, we estimate 

\begin{align*}
||u_{j}\cdot \nabla u_{j}||_{q}\leq ||u_j||_{\eta}||\nabla u_j||_{p}
\end{align*}\\
for $1/\eta=1/q-1/p$. Since $q\in [2,p]$ and $p>3$, we observe that $\sigma=3(1/p-1/\eta)=3(2/p-1/q)\leq 3/p<1$. By applying the interpolation inequality (B.2), we estimate

\begin{align*}
||u_j||_{L^{\eta}}
\leq C ||u_j||_{L^{p}}^{1-\sigma}||\nabla u_j||_{W^{1,p}}^{\sigma}.
\end{align*}\\
We take $T\leq 1$ and estimate 

\begin{align*}
||u_j\cdot \nabla u_{j}||_{L^{q}}
\leq C ||u_j||_{L^{p}}^{1-\sigma}||\nabla u_j||_{W^{1,p}}^{1+\sigma}
\leq C\Bigg(\frac{K_j}{s^{\gamma}}\Bigg)^{1-\sigma}\Bigg(\frac{2K_j}{s^{\gamma+\frac{1}{2}}}\Bigg)^{1+\sigma}
\leq C'\frac{K_j^{2}}{s^{\frac{3}{2}(1-\frac{1}{q}) }}.
\end{align*}\\
Since the above estimate holds for $s\leq T\leq 1$, we have 
\begin{align*}
||u_j\cdot \nabla u_j||_{\tilde{L}^{q}}
=\max\{||u_j\cdot \nabla u_j||_{L^{q}},||u_j\cdot \nabla u_j||_{L^{2}}   \}
\leq \frac{CK_j^{2}}{s^{\frac{3}{2}(1-\frac{1}{q}) }}.  \tag{A.4}
\end{align*} \\
Applying (A.1) implies 

\begin{align*}
||\partial_x^{k}e^{(t-s)A}\mathbb{P}u_j\cdot \nabla u_{j}||_{\tilde{L}^{p}}
\leq \frac{CK_j^{2}}{(t-s)^{\frac{3}{2}(\frac{1}{q}-\frac{1}{p})+\frac{|k|}{2}} s^{\frac{3}{2}(1-\frac{1}{q})} }.   \tag{A.5}
\end{align*}\\
We estimate $K_{j+1}$. We set $p_0=\max\{3p/(p+3), 2  \}$ and fix $q\in (p_0,3)$ so that the right hand-side of (A.5) is integrable near $s=t$ for $|k|\leq 1$. It follows from (A.1) and (A.5) that 

\begin{align*}
||u_{j+1}||_{\tilde{L}^{p}}
\leq ||e^{tA}u_0||_{\tilde{L}^{p}}
+\frac{C}{t^{\frac{1}{2}-\frac{3}{2p} }} K_{j}^{2}. 
\end{align*}\\
Similarly, we estimate $\nabla u_{j+1}$ and obtain $K_{j+1}\leq K_{1}+C_0 K_{j}^{2}$. Since the Stokes semigroup is strongly continuous on $\tilde{L}^{3}$, we have $K_{1}\to 0$ as $T\to0$. We take $T>0$ sufficiently small so that  $K_{1}\leq (4C_0)^{-1}$ and 

\begin{align*}
K_j\leq 2K_1\quad \textrm{for}\ j=1,2,\cdots.
\end{align*}\\
By a similar way, we estimate the difference $u_{j+1}-u_{j}$ and obtain 

\begin{align*}
\sup_{0\leq t\leq T}t^{\gamma}(||u_{j+1}-u_{j}||_{\tilde{L}^{p}}
+t^{\frac{1}{2}}||\nabla (u_{j+1}-u_{j})||_{\tilde{L}^{p}})\to 0\quad \textrm{as}\ j\to\infty.
\end{align*}\\
Thus the sequence $\{u_j\}$ converges to a mild solution $u$ satisfying (2.1) and (2.2) for $p, r\in (3,\infty)$. In particular, we have

\begin{align*}
K=\sup_{0\leq t\leq T}t^{\gamma}\{||u||_{\tilde{L}^{p}}+t^{\frac{1}{2}}||\nabla u||_{\tilde{L}^{p}}\}\leq 2K_1.   \tag{A.6}
\end{align*}\\
The uniqueness follows from the integral form since $t^{\gamma}u$ and $t^{\gamma+1/2}\nabla u$ vanish at time zero.\\

It remains to show (2.1) and (2.2) at end points. The property (2.1) for $p=\infty$ follows from the interpolation inequality (B.2). It follows from (A.1), (A.4) and (A.6) that 

\begin{align*}
||u-e^{tA}u_0||_{\tilde{L}^{3}}
&\leq \int_{0}^{t}||e^{(t-s)A}\mathbb{P}u\cdot \nabla u||_{\tilde{L}^{3}}\dd s\\
&\leq CK^{2}_{1}\int_{0}^{t}\frac{\dd s }{(t-s)^{\frac{3}{2}(\frac{1}{q}-\frac{1}{3}) }s^{\frac{3}{2}(1-\frac{1}{q}) }}
= C'K_1^{2}.
\end{align*}\\
Since $K_{1}\to 0$ as $T\to0$, the mild solution $u$ is strongly continuous on $\tilde{L}^{3}$ at time zero. Thus (2.1) holds for $p=3$. By a similar way, we estimate $t^{\frac{1}{2}}||\nabla u-\nabla e^{tA}u_0||_{\tilde{L}^{3}}\leq CK_1^{2}$. Since $t^{1/2}\nabla e^{tA}u_0$ vanishes on $\tilde{L}^{3}$ at time zero, (2.2) holds for $r=3$. 
\end{proof}

\vspace{5pt}

\begin{proof}[Proof of Lemma 2.1]
It remains to show the H\"older continuity (2.3). We set $f=-\mathbb{P}u\cdot \nabla u$. It follows from (A.4) that

\begin{align*}
||f||_{\tilde{L}^{q}}\leq \frac{CK^{2}}{s^{\frac{3}{2}(1-\frac{1}{q})}}\quad\textrm{for}\ 2\leq q\leq 3.\tag{A.7}
\end{align*}\\
We take an arbitrary $\delta \in (0,T)$ and $\alpha\in (0,1)$. For $\delta\leq \tau<t\leq T$, we estimate

\begin{align*}
||u(t)-u(\tau)||_{\tilde{L}^{3}}
&\leq ||e^{t A}u_0-e^{\tau A}u_0||_{\tilde{L}^{3}}
+\int_{\tau}^{t}||e^{(t-s)A}f||_{\tilde{L}^{3}}\dd s
+\int_{0}^{\tau}||e^{(t-s)A}f-e^{(\tau-s)A}f   ||_{\tilde{L}^{3}}\dd s \\
&=:I+II+III.
\end{align*}\\
It follows from (A.2), (A.1) and (A.7) that 

\begin{align*}
I&\leq C\delta^{-\alpha}(t-\tau)^{\alpha}||u_0||_{\tilde{L}^{3}},\\
II&\leq C\int_{\tau}^{t}||f||_{\tilde{L}^{3}}\dd s
\leq C'\delta^{-1}K^{2}(t-\tau).
\end{align*}\\
We estimate $III$. Since 

\begin{align*}
e^{(t-s)A}f-e^{(\tau-s)A}f=(e^{(t-\tau)A}-1)e^{\frac{(\tau-s)}{2}A}e^{\frac{(\tau-s)}{2}A}f,
\end{align*}\\
it follows from (A.2), (A.1) and (A.6) that 

\begin{align*}
|| e^{(t-s)A}f-e^{(\tau-s)A}f||_{\tilde{L}^{3}}
&\leq C\Big(\frac{t-\tau}{\tau-s}\Big)^{\alpha}||e^{\frac{(\tau-s)}{2}A}f ||_{\tilde{L}^{3}}\\
&\leq C' \frac{(t-\tau)^{\alpha}}{(\tau-s)^{\alpha+\frac{3}{2}(\frac{1}{q}-\frac{1}{3}) } }||f||_{\tilde{L}^{q}}\\
&\leq C''K^{2}\frac{(t-\tau)^{\alpha}}{(\tau-s)^{\alpha+\frac{3}{2}(\frac{1}{q}-\frac{1}{3}) }s^{\frac{3}{2}(1-\frac{1}{q}) }  }.
\end{align*}\\
We take $q\in [2,3)$ so that $3/2(1/q-1/3)<1-\alpha$ and obtain 

\begin{align*}
III\leq C\delta^{-\alpha}K^{2}(t-\tau)^{\alpha}.
\end{align*}\\
Thus $u\in C^{\alpha}([\delta,T]; \tilde{L}^{3})$ for $\alpha\in (0,1)$. By a similar way, $\nabla u\in C^{\alpha/2}([\delta,T]; L^{2})$ follows. We proved (2.3). The proof is now complete.
\end{proof}

\vspace{20pt}

%Appendix B
\section{Interpolation inequalities}

\vspace{10pt}

We give a proof for interpolation inequalities used in Proposition 5.2.

\vspace{10pt}

%lem B.1
\begin{lem}
The estimates 

\begin{align*}
||\varphi||_{p}&\leq C||\varphi||_{q}^{1-\sigma}||\nabla \varphi||_{q}^{\sigma},\quad \varphi\in W^{1,q}_{0},   \tag{B.1}  \\
||\phi||_{p}&\leq C||\phi||_{q}^{1-\sigma}|| \phi||_{1,q}^{\sigma},\quad \phi\in W^{1,q},    \tag{B.2}
\end{align*}\\
hold for $1\leq q\leq p\leq \infty$ satisfying $\sigma=3(1/q-1/p)<1$, where $W^{1,q}_{0}$ denotes the space of functions in $W^{1,q}$, vanishing on $\partial\Pi$.
\end{lem}

\vspace{5pt}

\begin{proof}
The estimate (B.1) for $\Pi=\mathbb{R}^{3}$ holds by estimates of the heat semigroup. Since the trace of $\varphi\in W^{1,q}_{0}$ vanishes on $\partial\Pi$, we apply (B.1) to the zero extension of $\varphi$ to $\mathbb{R}^{3}$ and obtain the desired estimate for $\Pi\subset \mathbb{R}^{3}$. For functions $\phi\in W^{1,q}$ with non-trivial traces, we use an extension operator $E: W^{1,q}(\Pi)\longrightarrow W^{1,q}(\mathbb{R}^{3})$ acting as a bounded operator also from $L^{q}(\Pi)$ to $L^{q}(\mathbb{R}^{3})$ \cite[Chapter VI, 3.1 Theorem 5]{Stein70}. By applying (B.1) for $\mathbb{R}^{3}$ and $E\phi$, we obtain (B.2).
\end{proof}

\vspace{15pt}

\section*{Acknowledgements}
The first author would like to thank Oxford University for their hospitality from October 2015 to January 2016. The first author was  supported by JSPS through the Grant-in-aid for Research Activity Start-up 15H06312, Young Scientist (B) 17K14217 and Scientific Research (B) 17H02853. The second author  was supported by the Ministry of Education and Science of the Russian Federation (grant 14.Z50.31.0037).

\vspace{15pt}

%ref
\bibliographystyle{plain}

\bibliography{ref}

\end{document}